\documentclass{article}
\def\WHO{nobody}
\def\version{25.3.2019}\def\users{}  %
\def\users{final-layout}  
\usepackage{amsmath,amsfonts,amssymb,color}
\usepackage{mathrsfs} 
\usepackage{epsfig}
\usepackage{psfrag}
\usepackage{upgreek} 
\usepackage{cite}
\usepackage{url}
\usepackage{import}
 
\newtheorem{theorem}{Theorem}[section]

\newtheorem{definition}[theorem]{Definition}

\newtheorem{proposition}[theorem]{Proposition}

\newtheorem{remark}[theorem]{Remark}

\numberwithin{equation}{section}
\usepackage{ifthen}
\ifthenelse{\equal{\users}{final-layout}}{}{
\usepackage{fancyhdr}
\pagestyle{fancy}
\headheight=28pt\headwidth=16.5cm
\definecolor{gray}{gray}{0.5}
\rhead{\color{gray}Explicit elastodynamics\\
T.Roub\'\i\v cek, C.G.Panagiotopoulos, C. Tsogka}
\chead{}
\lhead{Version \version, file: \jobname.tex, \underline{\Large\color{red}\WHO's working on} \\
compiled:
\number\day.\number\month.\number\year\ at
\the\hour:\ifnum\minute<10 0\fi\the\minute\ h }
}

\newcount\hour \newcount\minute
\hour=\time
\divide \hour by 60
\minute=\time
\loop \ifnum \minute > 59 \advance \minute by -60 \repeat

\usepackage[normalem]{ulem}
\definecolor{brown}{rgb}{0.5,0,0}
\ifthenelse{\equal{\users}{final-layout}}{
    \newcommand{\REPLACE}[2]{#2}
    
    \newcommand{\DELETE}[1]{}
    
    \newcommand{\COMMENT}[1]{}
    \newcommand{\COLOR}[1]{#1}
    \newcommand{\TINY}[1]{}
    \newcommand{\MARGINOTE}[1]{}
}
{\newcommand{\REPLACE}[2]{{\color{brown}\sout{#1}\uline{#2}\color{black}}}

 \newcommand{\DELETE}[1]{{\color{brown}\sout{#1}\color{black}}}
 
 \newcommand{\COMMENT}[1]{{\color{red}\uuline{#1}\color{black}}}
 \newcommand{\COLOR}[2][blue]{{\color{#1}{#2}}}
 \newcommand{\TINY}[1]{{\tiny#1}}
 \newcommand{\MARGINOTE}[1]{\marginpar{\color{red}\tiny\texttt{#1}}}
}

\renewcommand\dot[1]{\mathchoice
                 {{\buildrel{\hspace*{.1em}\text{\LARGE.}}\over{#1}}}
                 {{\buildrel{\hspace*{.1em}\text{\Large.}}\over{#1}}}
                 {{\buildrel{\hspace*{.1em}\text{\large.}}\over{#1}}}
                 {{\buildrel{\hspace*{.1em}\text{\large.}}\over{#1}}}}
\newcommand\DT{\dot}
\newcommand\DDT[1]{\mathchoice
   {{\buildrel{\hspace*{.1em}\text{\LARGE.\hspace*{-.1em}.}}\over{#1}}}
   {{\buildrel{\hspace*{.1em}\text{\Large.\hspace*{-.1em}.}}\over{#1}}}
   {{\buildrel{\hspace*{.1em}\text{\large.\hspace*{-.1em}.}}\over{#1}}}
   {{\buildrel{\hspace*{.1em}\text{\large.\hspace*{-.1em}.}}\over{#1}}}}

\newcommand{\wt}[1]{\mathchoice
     {\text{\small$\widetilde{\text{\normalsize$#1$}}\hspace*{.03em}$}}
                    {\text{\small$\widetilde{\text{\normalsize$#1$}}$}}
                    {\widetilde{#1\hspace*{-.02em}}\hspace*{.03em}}
                    {\widetilde{#1}}}

\def\R{{\mathbb R}}
\newcommand\calE{\mathcal E}
\newcommand\calH{\mathcal H}
\newcommand\calS{\mathcal S}
\newcommand\calU{\mathcal U}
\newcommand\calZ{\mathcal Z}
\newcommand\calX{\mathcal X}
\newcommand\bbC{\mathbb C}
\newcommand\bbD{\mathbb D}

\newcommand\bbB{\mathbb B}
\newcommand\bbI{\mathbb I}
\newcommand{\lineunder}[2]{\LU{\begin{array}[t]{c}\underbrace{#1}\vspace*{.5em}\end{array}}{\mbox{\footnotesize\rm #2}}}
\newcommand{\linesunder}[3]{\LSU{\begin{array}[t]{c}\underbrace{#1}\vspace*{.5em}\end{array}}{\mbox{\footnotesize\rm #2}}{\mbox{\footnotesize\rm #3}}}

\newcommand{\LU}[2]{\begin{array}[t]{c}#1\vspace*{-1em}\\_{#2}\end{array}}
\newcommand{\LSU}[3]{\begin{array}[t]{c}#1\vspace*{-1em}\\_{#2}\vspace*{-.5em}\\_{#3}\end{array}}

\newcommand{\nablaS}{\nabla_{\scriptscriptstyle\textrm{\hspace*{-.3em}S}}^{}}

\renewcommand\d{\mathrm d}

\newcommand{\KK}{k+1}
\newcommand{\KKK}{k}
\newcommand{\TT}{{\mathcal T}}

\newcommand{\barr}{\begin{array}}
\newcommand{\earr}{\end{array}}

\def\In{{\in}}
\newcommand{\eq}[1]{\eqref{#1}}

\begin{document}

\MARGINOTE{possible~journals:
\\
\\SIAM J.Sci.Comp.
\\or\\
SIAM Numer. Anal.
\\
or\\Numer. Math.
\\
or
\\.......???}
\begin{center}
\LARGE\bfseries
Explicit time-discretisation of elastodynamics
with some inelastic processes at small strains.
\end{center}

\bigskip

\begin{center}
\bfseries
Tom\'{a}\v{s} Roub\'\i\v{c}ek\footnote{Charles University,
Mathematical Institute, Sokolovsk\'a 83, CZ-186~75~Praha~8,
Czech Republic, {\tt tomas.roubicek@mff.cuni.cz}}${}^,$\footnote{Institute of
Thermomechanics of the Czech Acad. Sci., Dolej\v skova~5,
CZ--182~08 Praha 8, Czech Rep.},
Christos G. Panagiotopoulos\footnote{\label{fnote1}Institute of Applied and Computational Mathematics,
  Foundation for Research and Technology - Hellas, Nikolaou Plastira 100, Vassilika Vouton, GR-700 13
Heraklion, Crete, Greece
{\texttt{pchr@iacm.forth.gr}}},
Chrysoula Tsogka\footnote{Applied Math Unit, UC Merced, 5200 North Lake Road, Merced, CA 95343, USA, {\texttt{ctsogka@ucmerced.edu}}}${}^,$\footnotemark[3]
\end{center}

\bigskip
\noindent{\bf Abstract:} The 2-step staggered (also called leap-frog) time
discretisation of linear 2nd-order Hamiltonian systems (typically linear
elastodynamics in a stress-velocity form) is extended for a 3-step
staggered discretisation applicable for systems involving some internal
variables subjected to a dissipative evolution. After spatial discretisation,
a-priori estimates and convergence is proved under the usual CFL-condition.
Applications to specific problems in continuum mechanics of
solids at small stains are considered, in particular linearized plasticity, diffusion in
poroelastic media, damage, or adhesive contact. Numerical implementation
and some computational 2-dimensional simulation of waves emitted by a
rupture (delamination) of an adhesive contact illustrate the abstract
theory and efficiency of the explicit method.

\medskip


  \medskip


\section{Introduction -- mere linear elastodynamics}

In computational continuum mechanics of elastic or viscoelastic solids,
so-called transient 
(low-frequency) and wave propagation problems(high-frequency) are 
distinguished and different approximation methods are used.
The prototype equation or rather initial-boundary-value problem which we have
in mind is the linear elastodynamic at small strains:
\begin{subequations}\label{IBVP}\begin{align}\label{IBVP1}
    &  \varrho\DDT u-{\rm div}\,\bbC e(u)=f&&\text{on }\ \varOmega\
    \text{ for }\ t\in[0,T],
  \\&\label{IBVP2}
    [\bbC e(u)]\vec{n}+\bbB u=g&&\text{on }\ \varGamma\
    \text{ for }\ t\in[0,T],
  \\&\label{IBVP3}
  u|_{t=0}=u_0,\ \ \ \DT u|_{t=0}=v_0&&\text{on }\ \varOmega,
\end{align}\end{subequations}
with $\varrho>0$ a mass density, $\bbC$ a symmetric positive definite
4th-order elasticity tensor, $e(u)=\frac12(\nabla u)^\top{+}\frac12\nabla u$
the small-strain tensor, $\bbB$ a symmetric positive semidefinite
2nd-order tensor determining the elastic support on the boundary,
$f$ the bulk force, $g$ the surface loading,
$u_0$ the prescribed initial
displacement, $v_0$ initial velocity, and $T>0$ a fixed time horizon.
The unknown $u:[0,T]\time\varOmega\to\R^d$ is the displacement, the
dot-notation stands for the time derivative,
$\Omega\subset\R^d$ is a bounded Lipschitz domain, $d=2$ or $3$,
$\varGamma$ is
its boundary, and $\vec{n}$ the unit outward normal.
For notational simplicity, we will write the initial-boundary-value
problem \eqref{IBVP} in the abstract form 
\begin{align}\label{IVP}
  \mathcal T'\DDT u+\mathcal W' u=\mathcal F_u'(t)\ \ 
\text{ for }\ t\in[0,T],\ \ u|_{t=0}=u_0,\ \ \ \DT u|_{t=0}=v_0,
\end{align}
where $\mathcal T$ is the kinetic energy, $\mathcal W$ is the stored
energy, and $\mathcal F$ is the external force, while
$(\cdot)'$ denotes the G\^ateaux derivative. In the context
of \eqref{IBVP}, $\mathcal T(v)=\int_\varOmega\frac12\varrho|v|^2\,\d x$,
$\mathcal W(u)=\int_\varOmega\frac12\bbC e(u){:}e(u)\,\d x
+\int_\varGamma\frac12\bbB u\cdot u\,\d S$, and
$\mathcal F(t,u)=\int_\varOmega f(t){\cdot}u\,\d x
+\int_\varGamma g(t){\cdot}u \,\d S$. Thus
$\mathcal F_u'(t)$ is the linear functional, let us denote it shortly by
$F(t)$.

In situations where high-frequency oscillations arising
  typically during wave propagation are to be calculated,
the implicit time-discretisations (even if energy conserving as
e.g.\ \cite{PaPaTa09RDCT,RouPan17ECTD})
are computationally cumbersome especially in 3-dimensional 
problems. Hence \emph{explicit time discretisations} 
are more efficient. The simplest explicit scheme is the so-called
  central-difference scheme
\begin{align}\label{elast-dyn-explicit}
\mathcal T'\frac{u_{\tau h}^{k+1}-2u_{\tau h}^{k}+u_{\tau h}^{k-1}}{\tau^2}
+\mathcal W_h' u_{\tau h}^{k}=F_h(k\tau)
\end{align}
with $\tau>0$ a time step, and with
$\mathcal W_h$ and
$F_h$ denoting some approximations of the respective functionals
obtained by a suitable finite-element method (FEM) with the mesh size $h>0$.
In particular, a numerical approximation 
leading to a diagonalization of the mass matrix $\mathcal T'$, called mass lumping,  
in \eq{elast-dyn-explicit} is an important ingredient so as to obtain 
efficient explicit methods.
The formula \eqref{elast-dyn-explicit} leads, when tested by 
$\frac{u_{\tau h}^{k+1}-u_{\tau h}^{k}}{\tau}$, to a correct discrete kinetic 
energy $\mathcal T$ but a twisted stored energy,
namely
$\frac12\langle\mathcal W_h'u_{\tau h}^{k+1},u_{\tau h}^{k}\rangle$,
whose handling needs
the \emph{Courant-Fridrichs-Lewy (CFL) condition} \cite{CoFrLe28PDMP} that
typically bounds the time discretisation step $\tau=\mathscr{O}(h_{\rm min})$ with $h_{\rm min}$ the smallest
element size on a FEM discretisation 
cf.\ \cite[p.171]{Tsog99} or also e.g.\ \cite{KPCOP16TSDS,NPSK07SDAE}.
More specifically, the CFL reads as
\begin{align}
  \big\langle \mathcal T' u_h,u_h\big\rangle 
  \geq \frac{\tau^2}{4}\big\langle{\cal W}_h' u_h,u_h \big\rangle
  \label{CFLold}
\end{align}
for any $u_h$ from the respective finite-dimensional subspace.
This is a drawback which makes such discretisation
less suitable for enhancing the stored energy by some internal variables and
(possibly) nonlinear processes on them, which is the goal of this article.

Therefore, we use another, so-called \emph{leap-frog}, \emph{scheme}. To this
aim, we first rewrite \eqref{IBVP1} in in the velocity/stress formulation,
i.e.\ terms of $v=\DT u$ and of
the stress $\sigma:=\bbC e(u)$, eliminating the displacement $u$, further we consider the rate form of \eqref{IBVP2}, together with appropriate initial conditions:
\begin{subequations}\label{IBVP+}\begin{align}\label{IBVP1+}
    &&&  \varrho\DT v-{\rm div}\,\sigma=f\ \ \ \text{ and }\ \ \ 
    \DT\sigma=\bbC e(v)&&\text{on }\ \varOmega\ \text{ for }\
  t\in[0,T],&&
  \\&&&\label{IBVP2+}
  \DT\sigma\vec{n}+\bbB 
  v
  =\DT g&&\text{on }\ \varGamma\ \text{ for }\
  t\in[0,T],
  \\&&&\label{IBVP3+}
    v|_{t=0}=v_0,\ \ \ \sigma|_{t=0}=\sigma_0:=\bbC e(u_0)&&\text{on }\ \varOmega.
\end{align}\end{subequations}
In the abstract form \eqref{IVP}, when writing $\mathcal W=\mathscr{W}\circ
E$ with $E$ denoting the linear operator $u\mapsto(e,w):=(e(u),u|_\varGamma)$,
this reads as 
\begin{subequations}\label{IVP+}\begin{align}
  &&&\mathcal T'\DT v+E^*\varSigma=F(t)&&\text{ for }\ t\in[0,T],&&
    v|_{t=0}^{}=v_0,\ \ 
  \ \ \text{ and }\ \
  \\&&&\DT\varSigma=\mathscr{W}'Ev+\DT G(t)  
  &&\text{ for }\ t\in[0,T],
    &&\varSigma|_{t=0}^{}=\varSigma_0:=\mathscr{W}'Eu_0,&&&&\label{IBVP2++}
\end{align}\end{subequations}
where $E^*$ is the adjoint operator to $E$. The stored energy
governing \eqref{IBVP+} is
$\mathscr{W}(e,w)=\int_\varOmega\frac12\bbC e{:}e\,\d x+
\int_\varGamma\frac12\bbB w{\cdot}w\,\d S$ while
the external loading is now split into two parts acting differently,
namely $\langle F(t),u\rangle=\int_\Omega f(t)\cdot u\,\d x$ and
$\langle G(t),w\rangle=\int_\varGamma g(t)\cdot w\,\d x.$ 
Let us note that \eqref{IVP+} involves, in fact, the
equation on $\varOmega$ as well as the equation on $\varGamma$ if
$\mathcal T$ understood as the functional on $\varOmega\times\varGamma$,
being trivial on $\varGamma$ since no inertial is considered on
the $(d{-}1$)-dimensional boundary $\varGamma$. 
In particular, the ``generalized'' stress
$\varSigma=\mathscr{W}'Eu=(\bbC e(u),\bbB u|_\varGamma)$
contains, beside the bulk stress tensor, also the traction
stress vector. Relying on the linearity of $\mathscr{W}'$,
we have $\DT\varSigma=\mathscr{W}'Ev$ with $v=\DT u$, as used in
\eqref{IVP+}.

The mentioned ``leap-frog'' time discretisation of
\eqref{IVP+} then reads as 
\begin{align}\label{elast-dyn-explicit+}
  \frac{\varSigma_{\tau h}^{k+1/2}-\varSigma_{\tau h}^{k-1/2}}\tau=\mathscr{W}'E_hv_{\tau h}^{k}+D_{\tau h}^{k}
  \ \ \ \ \ \text{ and }\ \ \ \
  \TT'\frac{v_{\tau h}^{k+1}{-}v_{\tau h}^k}\tau+E_h^*\varSigma_{\tau h}^{k+1/2}
  =F_{\tau h}^{k+1/2},
\end{align}
where $\mathscr{W}_h$
and $E_h$ is a suitable FEM discretisation of $E$
and
\begin{align}\label{F-G}
  F_{\tau h}^{k+1/2}:=\frac1\tau\int_{k\tau}^{(k+1)\tau}\!\!\!\!\!\!\!\!\!\!\!F_h(t)\,\d t\
  \text{ and }\ \
  D_{\tau h}^{k}:=\frac1\tau\int_{(k-1/2)\tau}^{(k+1/2)\tau}\!\!\DT G_h(t)\,\d t
  =\frac{G_h((k{+}\frac12)\tau)-G_h((k{-}\frac12)\tau)}\tau.
\end{align}
We assume that $v$'s is discretised in a piecewise-constant way so that
$\mathcal T$ leads to a diagonal form on such a subspace and therefore numerical
integration leading to mass lumping is not needed here. Otherwise,
higher-order discretisation with mass lumping may also be used to achieve the
desired property of obtaining an explicit scheme (avoid solving systems of equations). 
We refer to \cite{BeJoTs02NFMF,JolTso08FEMD,Tsog99} for details in the case $G\equiv0$. 
This discretisation
also does not need inversion of $\mathcal W_h'=E_h^*\mathscr{W}'E_h^{}$,
which is just the ultimate goal of all explicit discretisation schemes.
Usually, the spatial FEM discretisation exploits regularity available in
linear elastodynamics, in particular that ${\rm div}\,\sigma$ and
$e(v)$ in \eqref{IBVP1+} live in $L^2$-spaces.
Moreover, the equations in \eqref{elast-dyn-explicit+} are decoupled in the
sense that, first, $\varSigma_{\tau h}^{k+1/2}$ is calculated from
the former equation and, second, $v_{\tau h}^{k+1}$ is calculated from the latter
equation assuming, that $(v_{\tau h}^k,\varSigma_{\tau h}^{k-1/2})$ is known
from the previous time step. For $k=0$, it starts from $v_{\tau h}^0=v_0$
and from a  half time step  
$\varSigma_{\tau h}^{1/2}=\varSigma_{\tau h}^{0}+\frac{\tau}{2} \mathscr{W}'E_hv_{\tau h}^{0}$.
For the space discretisation, the lower order $Q^{\rm div}_{k+1}-Q_{k}$ finite element is obtained for $k=0$ and in this case the velocity is discretised as piecewise
constant on rectangular or cubic elements while the stress is discretised by piecewise bi-linear functions with
some continuities. Namely the normal component of the stress is continuous
across edges of adjacent elements while the tangential component is allowed
to be discontinuous. For more details about the space discretisation we
refer the interested reader to \cite{BeJoTs02NFMF}. 
An alternative discretisation using triangular elements known as
{\it staggered discontinuous Galerkin method} is proposed in
\cite{ChLaQi15SDGM}. In general, the leap-frog scheme has been frequently used in
  geophysics to calculate seismic wave propagation with the finite differences method, cf.\ e.g.\
\cite{Bohl02PVFD,Grav96SSWP,Viri84SHWP}. 

When taking the average (i.e.\ the sum with the weights $\frac12$ and
$\frac12$) of the second equation in \eqref{elast-dyn-explicit+}
in the level $k$ and $k{-}1$ tested by $v_{\tau h}^k$ and summing it with the
first equation in \eqref{elast-dyn-explicit+} tested by
$[\mathscr{W}']^{-1}(\varSigma_{\tau h}^{k+1/2}{+}\varSigma_{\tau h}^{k-1/2})/2$,
we obtain
\begin{align*}&\frac1{2\tau}\big\langle[\mathscr{W}']^{-1}\varSigma_{\tau h}^{k+1/2},
\varSigma_{\tau h}^{k+1/2}\big\rangle
-\frac1{2\tau}\big\langle[\mathscr{W}']^{-1}\varSigma_{\tau h}^{k-1/2},
\varSigma_{\tau h}^{k-1/2}\big\rangle
\\&\qquad\qquad
=\Big\langle\frac{\varSigma_{\tau h}^{k+1/2}{+}\varSigma_{\tau h}^{k-1/2}}2,
E_h^{}v_{\tau h}^{k}\Big\rangle+\Big\langle
[\mathscr{W}']^{-1}D_{\tau h}^k,
\frac{\varSigma_{\tau h}^{k+1/2}{+}\varSigma_{\tau h}^{k-1/2}}2\Big\rangle
\ \ \text{ and}
\\
&\Big\langle\TT'\frac{v_{\tau h}^{k+1}{-}v_{\tau h}^{k-1}}\tau,
v_{\tau h}^{k}\Big\rangle
+\Big\langle\frac{\varSigma_{\tau h}^{k+1/2}{+}\varSigma_{\tau h}^{k-1/2}}2,E_h^{}v_{\tau h}^{k}\Big\rangle
=\Big\langle\frac{ F_{\tau h}^{k+1/2}{+}F_{\tau h}^{k-1/2}}2,v_{\tau h}^{k}\Big\rangle\,.
\end{align*}
Summing it up, we eventually obtain the (approximate) energy balance with
the correct  stored energy and twisted discrete kinetic energy,
namely 
\begin{align}\label{leap-frog-energy}
\frac12\langle\mathscr{T}'v_{\tau h}^{k+1},v_{\tau h}^k\rangle+
\varPhi_h
(\varSigma_{\tau h}^{k+1/2})\ \ \text{ with }\ \ 
\varPhi_h(\varSigma)=\frac12\langle[\mathscr{W}']^{-1}\varSigma,\varSigma\rangle\,;
\end{align}
note that $\varPhi_h$ is the (possibly approximate) stored energy but expressed in terms of the
generalized stress. Yet, in contrast to  \eqref{elast-dyn-explicit}, 
yielding the energy balance with the correct stored energy,
\eqref{elast-dyn-explicit+} allows for enhancement of this stored energy by
some internal variables. This last attribute is a qualitative difference
compared to \eqref{elast-dyn-explicit}.
Again, the a-priori estimates
and convergence for $\tau\to0$ and $h\to0$ needs the following CFL condition
\begin{align}
  \big\langle[ \mathcal W'_h ]^{-1}\varSigma_h,\varSigma_h\big\rangle
  \geq \frac{\tau^2}{4}\big\langle E^*_h \varSigma_h,(\TT')^{-1}  E^*_h \varSigma_h\big\rangle 
  \label{CFLold2}
\end{align}
for any $\varSigma_h$ from the respective finite-dimensional subspace.
Moreover, $F=0$ is often considered, which makes
the a-priori estimation easier.
Let us also note that the adjective ``leap-frog'' is sometimes
used also for the time-discretisation \eqref{elast-dyn-explicit} if
written as a two-step scheme, cf.\ e.g.\ \cite[Sect.\,7.1.1.1]{CohPer17FEDG}.

The plan of this article is as follows: In Section~\ref{sec-int-var},
we extend the abstract system \eqref{IVP+} by another equation for
some internal variable and cast its weak formulation without relying on any regularity.
Then, in Section~\ref{sect-disc}, we enhance the
two-step leap-frog discrete scheme \eqref{elast-dyn-explicit+}
to a suitable three-step scheme, and show its energetics.
Then, in Section~\ref{sect-anal}, we   
prove the numerical stability of the 3-step staggered approximation scheme
and its convergence under the CFL condition modified correspondingly,
cf.\ \eqref{CFLnew} below.
Such an abstract scheme is then illustrated in Section~\ref{sec-exa}
on several examples from continuum mechanics, in particular
on models of plasticity, creep, diffusion, damage, and delamination.
Eventually, in Section~\eqref{sec-numer},  numerical implementation of the presented scheme for problems of adhesive contact is considered and computational experiments are shown in order to demonstrate its computational efficiency.

It should be emphasized that, to the best of our knowledge, a rigorously
justified (as far as numerical stability and convergence) combination of the
explicit staggered discretisation with nonlinear dissipative processes
on some internal variables is new, 
although occasionally some dissipative nonlinear phenomena can be found in
literature as in \cite{Scar04ETNP} for a unilateral contact, 
  in \cite{Bohl02PVFD} for a Maxwell viscoelastic rheology,
in \cite{SePaPa18SEIF} for electroactive polymers,
or in \cite{FePaFa01PACM} for general thermomechanical systems,
but without any numerical stability (a-priori estimates) and
convergence guaranteed. 

\COMMENT{MAYBE ALSO:
  T.J.R. Hughes, W.-K. Liu, Implicit-explicit finite elements in transient
  analysis: I. Stability theory; II. Implementation andnumerical examples,
  J. Appl. Mech. 45 (1978) 371--378.
  \\OR\\
  T.J.R. Hughes, K.S. Pister, R.L. Taylor, Implicit-explicit
  finite elements in nonlinear transient analysis, Comput. Methods Appl.Mech. Engrg. 17/18 (1979) 159--182
  \\OR\\
  T.J.R. Hughes, R.S. Stephenson, Stability of implicit-explicit finite elements in nonlinear transient analysis, Int. J. Engrg. Sci. 199 (1981) 295--302}

\COMMENT{Tomas if you want you can add the last of these references.} 

\section{Internal variables and their dissipative evolution.}\label{sec-int-var}

The concept of internal variables has a long tradition
and opens wide options for material modelling while the internal
parameters are subjected to 1st-order evolution flow rules,
cf.\ \cite{Maug15SIVS}. The system \eqref{IVP} is thus enhanced as:
\begin{subequations}\label{IVP2}
  \begin{align}\label{IVP2Sigma}
  &\mathcal T'\DDT u+\mathcal W_u'(u,z)=F(t)&&
  \text{ for }\ t\in[0,T],\ \ u|_{t=0}=u_0,\ \ \ \DT u|_{t=0}=v_0,
  \\&\label{IVP2z}\partial\varPsi(\DT z)+\mathcal W_z'(u,z)\ni0&&
  \text{ for }\ t\in[0,T],\ \ z|_{t=0}=z_0.
\end{align}\end{subequations}
The inclusion in \eqref{IVP2z} refers to a possibility that the convex
(pseudo)potential of dissipative forces $\varPsi$ may be nonsmooth and then its
subdifferential $\partial\varPsi$ can be multivalued.

Combination of the 2nd-order evolution
\eqref{IVP} with such 1st-order evolution is to be made carefully.
In contrast to the implicit schemes, cf.\ \cite{RouPan17ECTD},
the constitutive equation is differentiated in time, cf.\ \eqref{IBVP1+},
and it seems necessary to use the split (staggered) scheme so that
the internal-variable flow rule can be used without being
differentiated in time, even if the stored energy $\mathcal W$ would
be quadratic.

Moreover, to imitate the leap-frog scheme, it seems suitable (or maybe
even necessary) that the stored energy $\mathcal W$ can be expressed
in terms of the generalized stress as
\begin{align}\label{ansatz}
  \mathcal W(u,z)=\varPhi(\varSigma,z)\ \ \text{ with }\ \ 
  \varSigma=\mathfrak C Eu,\ \text{ and }\ \ \varPhi(\cdot,z)\ \text{ and }\
  \varPhi(\varSigma,\cdot)\text{ quadratic},
\end{align}
where $\mathfrak C$ stands for a ``generalized'' elasticity
tensor and $E$ is an abstract gradient-type operator;
  typically $Eu=(e(u),u|_\varGamma)$ or also simply $Eu=e(u)$ are
here considered in the context of continuum mechanics at small strains,
cf.\ the examples in Sect.\,\ref{sec-exa}.
Here, $\varSigma$ may not directly enter the balance of forces and
is thus to be called rather as some ``proto-stress'',
while the actual generalized stress will be denoted by $S$.
For a relaxation of the last requirement of \eqref{ansatz} see Remark~\ref{rem-nonquadratic}
below.

Then, likewise \eqref{IVP+}, we can write the system \eqref{IVP2}
in the velocity/proto-stress formulation as
\begin{subequations}\label{IVP2+}
  \begin{align}\label{IVP2Sigma+}
&\DT\varSigma=\mathfrak CEv+\DT G(t)
        &&
    \text{ for }\ t\in[0,T],\ \ &&
    \\ \label{IVP2S+}   &\mathcal T'\DT v+E^*S 
    =F(t)\ \ \text{ with }\ \
    S=\mathfrak C^*\varPhi_\varSigma'(\varSigma,z)\!\!\!\!\!\!&&
    \text{ for }\ t\in[0,T],\ \ 
    \\&\label{IVP2z+}
    \partial\varPsi(\DT z)+\varPhi_z'(\varSigma,z)  
    \ni0&&
    \text{ for }\ t\in[0,T],
    \\&\label{IVP2+IC}
    \varSigma|_{t=0}=\varSigma_0:=\mathfrak CEu_0+G(0)\,,\ \ \ v|_{t=0}=v_0\,,
    \ \ z|_{t=0}=z_0\,.\hspace*{-9em}
\end{align}\end{subequations}
Here $\varPhi_\varSigma'(\varSigma,z)$ is in a position of a ``generalized''
strain and, when multiplied by $\mathfrak C^*$, it becomes a 
generalized stress.

The energetics of this system can be revealed by testing the particular
equations/inclusions in \eqref{IVP2+} by
$\varPhi_\varSigma'(\varSigma,z)$, $v$, and $\DT z$.
Thus, at least formally, we obtain
\begin{subequations}\label{energy-test}\begin{align}
&
 \big\langle\varPhi_\varSigma'(\varSigma,z),\DT\varSigma
\big\rangle=\big\langle\varPhi_\varSigma'(\varSigma,z),\mathfrak CEv+\DT G\big\rangle
=\big\langle\mathfrak C^*\varPhi_\varSigma'(\varSigma,z),Ev\big\rangle
+\big\langle\varPhi_\varSigma'(\varSigma,z),\DT G\big\rangle,
  \\  &\mathcal T(v)+\big\langle\mathfrak C^*\varPhi_\varSigma'(\varSigma,z),Ev\big\rangle
  =\big\langle F(t),v\big\rangle\,,
  \\&\varXi(\DT z)
    +\big\langle\varPhi_z'(\varSigma,z),\DT z\big\rangle\le0\ \ \ \
  \text{ with }\ \ \ 
\varXi(\DT z):=\inf\big\langle\partial\varPsi(\DT z),\DT z\big\rangle\,.
\label{energy-test3}
\end{align}\end{subequations}
The functional $\varXi$ is in the position of the dissipative rate and 
the ``inf'' in it refers to the fact that the dissipative
potential $\varPsi$ can be nonsmooth and thus the 
subdifferential $\partial\varPsi$ can be multivalued even at $\DT z\ne0$,
otherwise an equality in \eqref{energy-test3} holds.
Summing it up and using
the calculus $\frac{\d}{\d t}\varPhi(\varSigma,z)=
\langle\varPhi_\varSigma'(\varSigma,z),\DT\varSigma\rangle
+\langle\varPhi_z'(\varSigma,z),\DT z\rangle$,
we obtain the energy (im)balance
\begin{align}\label{energy}
  \frac{\d}{\d t}\big(\!\!\!\linesunder{\mathcal T(v)
    +\varPhi(\varSigma,z)}{kinetic and}{stored energies}\!\!\!\big)\
  +\!\!\!\!\linesunder{\varXi(\DT z)
      }{dissipation}{rate}\!\!\!\!
  =\!\!\linesunder{\big\langle F(t),v\big\rangle
   +\big\langle\varPhi_\varSigma'(\varSigma,z),\DT G\big\rangle}
  {power of}{external force}\,.
\end{align}

Let us now formulate some abstract functional setting of the system
\eqref{IVP2+}. For some Banach spaces $\calS$, $\calZ$, and
$\calZ_1\supset\calZ$
and for a Hilbert space $\calH$, let 
$\varPhi:{\calS}\times\calZ\to\R$ be smooth and coercive, $\mathcal T:\calH\to\R$
be quadratic and coercive, and let $\Psi:\calZ\to[0,+\infty]$ 
be convex, lower semicontinuous, and coercive on $\calZ_1$, cf.\ \eqref{coerc+Lipschitz} below.
Intentionally, we do not want to rely on any regularity which is
usually at disposal in linear problems but might
be restrictive in some nonlinear problems. For this reason, we reconstruct
the abstract ``displacement'' and use \eqref{IVP2Sigma+} integrated in time,
i.e.
\begin{align}\label{IVP2Sigma++}
  \varSigma={\mathfrak C}Eu+G\ \ \ \text{ with }\ \
  u(t):=\int_0^t\!\!v(t)\,\d t+u_0\,.
\end{align}
Moreover, we still need another Banach space $\calE$
and define the Banach space $\calU:=\{u\in \calH;\
Eu\in \calE\}$ equipped with the standard graph norm. Then, by definition,
we have the continuous embedding $\calU\to\calH$ and the
continuous linear operator $E:\calU\to\calE$. We assume that
$\calU$ is embedded into $\calH$ densely, so that
$\calH^*\subset\calU^*$ and
that $\calH$ is identified with its dual $\calH^*$, so that
we have the so-called Gelfand triple
$$
\calU\subset\calH\cong\calH^*\subset\calU^*.
$$
We further consider the abstract elasticity tensor $\mathfrak C$
as a linear continuous operator $\calE\to\calS$.
Therefore $\mathfrak C Eu\in\calS$ provided $u\in \calU$ so that
the equation \eqref{IVP2Sigma++} is meant in $\calS$ and one needs
$G(t)\in\calS$.
Let us note that $\mathcal T':\calH\to\calH^*\cong\calH$,
$\varPhi_\varSigma':{\calS}\times\calZ\to{\calS}^*$,
$E^*:\calE^*\to\calU^*$, and
$\mathfrak C^*:\calS^*\to\calE^*$, so that
 $\mathcal T' v\in\calH^*$
provided $v\in\calH$
and also $S=\mathfrak C^*\varPhi_\varSigma'\in\calE^*$ and
$E^*S\in\calH^*$. In particular, the equation
\eqref{IVP2S+} can be meant in $\calH$ if integrated in time,
and one needs $F(t)$ valued in $\calH$.

We will use the standard notation $L^p(I;\calX)$ for Bochner spaces of
Bochner measurable functions $I\to\mathcal X$ whose norm is integrable with
the power $p$ or essentially bounded if $p=\infty$, and $W^{1,p}(I;\calX)$
the space of functions from $L^p(I;\calX)$ whose distributional time derivative
is also in  $L^p(I;\calX)$. Also, $C^k(I;\calX)$ will denote the space of
functions $I\to\calX$ whose $k^{\rm th}$-derivative is continuous, and
$C_{\rm w}(I;\calX)$ will denote the space of weakly continuous functions
$I\to\mathcal X$. Later, we will also use ${\rm Lin}(\calU,\calE)$,
denoting the space of linear bounded operators $\calU\to\calE$
normed by the usual sup-norm.

A weak formulation of \eqref{IVP2S+} can be obtained after by-part integration
over the time interval $I=[0,T]$ when tested by a smooth function.
It is often useful to confine ourselves to situations
\begin{align}\label{z-ansatz}
  \varPhi(\varSigma,z)=\varPhi_0(\varSigma,z)+\varPhi_1(z)
  \ \ \ \text{ with }\ [\varPhi_0]_z':{\calS}\times\calZ
  \to\calZ_1^*\ \text{ and }\ \varPhi_1':\calZ\to\calZ^*
\end{align}
and to use a by-part integration for the term
$\langle\varPhi_1'(z),\DT z\rangle$. Altogether, we arrive to:

\begin{definition}[Weak solution to \eqref{IVP2+}.]\label{def}
  The quadruple $(u,\varSigma,v,z)\in C_{\rm w}(I;{\calU})\times
  C_{\rm w}(I;{\calS})\times C_{\rm w}(I;{\calH})\times C_{\rm w}(I;\calZ)$
with $\varPsi(\DT z)\in L^1(I)$ 
  will be called a weak solution to the initial-value problem
  \eqref{IVP2+} with \eqref{IVP2Sigma++}
  if $v=\DT u$ in the distributional sense, $\varSigma={\mathfrak C}Eu+G$
  holds a.e.\ on $I$, and if
 \begin{subequations}\label{weak}\begin{align}
&\label{weak-v}
\int_0^T\!\!\big\langle \varPhi_\varSigma'(\varSigma,z),\mathfrak C E\wt v\big\rangle_{{\calS}^*\times{\calS}}^{}
    -\big\langle\mathcal T'v,\DT{\wt v}\big\rangle_{{\calH}^*\times{\calH}}^{}
    \,\d t=\big\langle\mathcal T'v_0,\wt v(0)\big\rangle_{{\calH}^*\times{\calH}}^{}
    +\int_0^T\!\!\big\langle F,\wt v\big\rangle_{{\calH}^*\times{\calH}}^{}\,\d t
    \intertext{for any $\wt v\in C^1(I;{\calH})\,\cap\, C(I;{\calU})$ with $\wt v(T)=0$, and}
   &\nonumber
   \int_0^T\!\!\varPsi(\wt z)+\big\langle[\varPhi_0]_z'(\varSigma,z),\wt z{-}
    \DT z\big\rangle_{\calZ_1^*\times\calZ_1}^{}\!
    +\big\langle\varPhi_1'(z),\wt z\big\rangle_{\calZ^*\times\calZ}^{}
        \,\d t+\varPhi_1(z_0)
       \\[-.6em]&\hspace*{22em}
    \ge\varPhi_1(z(T))+
       \!\int_0^T\!\!\varPsi(\DT z)\,\d t
    \label{weak-z}
 \end{align}\end{subequations}
 for any $\wt z\in C(I;\calZ)$, where indices in the dualities
 $\langle\cdot,\cdot\rangle$ indicate the respective spaces in dualities, and
 if also $u(0)=u_0$, $\varSigma(0)=\varSigma_0$, and $z(0)=z_0$.
\end{definition}

Let us note that the remaining initial condition $v(0)=v_0$ is contained in
\eqref{weak-v}. This definition works successfully for $p>1$, i.e.\ for
rate-dependent evolution of the abstract internal variable $z$, so that
$\DT z\in L^p(I;\calZ_1)$. For the rate-dependent evolution when $p=1$,
we would need to modify it but we will need to restrict ourselves for $p\ge2$,
see due to the a-priori estimates in Proposition~\ref{prop1}.

\section{A three-step staggered time discretisation}\label{sect-disc}

Now we devise the leap-frog discretisation of (\ref{IVP2+}a,b) combined with
the fractional-step split (a staggered scheme) with a mid-point
formula for \eqref{IVP2z+}. Instead of a two-step formula
\eqref{elast-dyn-explicit+}, we will obtain a three-step formula and
therefore, from now on, we will leave the convention
of a half-step notation used standardly in \eqref{elast-dyn-explicit+}
and write $k+1$ instead of $k+1/2$. 
Considering that we know from previous step
$\varSigma_{\tau h}^{\KKK},v_{\tau h}^k,z_{\tau h}^{k}$, then it leads to:
\begin{subequations}\label{suggestion}
\begin{align}\label{suggestion-1}
&\text{1) calculate $\varSigma_{\tau h}^{\KK}$:}&&\!\!\!
\frac{\varSigma_{\tau h}^{\KK}-\varSigma_{\tau h}^{\KKK}}\tau=
\mathfrak CE_h v_{\tau h}^{k}+D_{\tau h}^{k},
\\\label{suggestion-3}&\text{2) calculate $z_{\tau h}^{k+1}$:}&&\!\!\!
\partial\varPsi\Big(\frac{z_{\tau h}^{k+1}{-}z_{\tau h}^{k}}\tau\Big)+
\varPhi_z'\Big(\varSigma_{\tau h}^{\KK},\frac{z_{\tau h}^{k+1}{+}z_{\tau h}^{k}}{2}\Big)\ni0\,,
\\\nonumber
  &\text{3) calculate $v_{\tau h}^{k+1}$:}&&\!\!\! 
  \TT'\frac{v_{\tau h}^{k+1}{-}v_{\tau h}^{k}}{\tau}+E_h^\ast S_{\tau h}^{\KK}
  =F_{\tau h}^{\KK}
  \ \text{ with }\ S_{\tau h}^{\KK}=\mathfrak C^*\varPhi_\varSigma'(\varSigma_{\tau h}^{\KK},\COLOR{z_{\tau h}^{k+1}}),&&
\label{suggestion-2}
\\&\text{\ \ \ \ \ \ \ \ \ and \ $u_{\tau h}^{k+1}$:}&&\!\!\!
u_{\tau h}^{k+1}=u_{\tau h}^{k}+\tau v_{\tau h}^{k+1},
\end{align}\end{subequations}
where $F_{\tau h}^{\KK}$ and $D_{\tau h}^{k}$ are from \eqref{F-G}. 

The only possibly nonlocal equation is \eqref{suggestion-3}. This equation
has a potential
\begin{align}
z\mapsto\frac2\tau\varPhi\Big(\varSigma_{\tau h}^{\KK},\frac{z{+}z_{\tau h}^{k}}{2}\Big)+\varPsi\Big(\frac{z{-}z_{\tau h}^{k}}\tau\Big)
\end{align}
and therefore the existence of a solution to the
inclusion or rather variational inequality
\eqref{suggestion-3} can be shown by a direct method,
cf.\ also \cite{RouPan17ECTD}. In view of the definition of the
convex subdifferential, \eqref{suggestion-3} means the variational
inequality
\begin{align}\label{VI}
  \varPsi\big(\widetilde z\big)+\bigg\langle
\varPhi_z'\Big(\varSigma_{\tau h}^{\KK},\frac{z_{\tau h}^{k+1}{+}z_{\tau h}^{k}}{2}\Big),
  \widetilde z-\frac{z_{\tau h}^{k+1}{-}z_{\tau h}^{k}}\tau\bigg\rangle
  \ge\varPsi\Big(\frac{z_{\tau h}^{k+1}{-}z_{\tau h}^{k}}\tau\Big)
\end{align}
for any $\widetilde z$.

To allow a discontinuous Galerkin discretisation for the
velocity like in \cite{BeJoTs02NFMF,BeJoTs01MFES,ChLaQi15SDGM,Tsog99},
we consider $E_h\in{\rm Lin}(\calH,\calE)$. Later, we will need
the approximation property 
\begin{align}\label{Eh->E}
  v\in\calU,\ v_h\in V_h,\ \ v_h\to v\ \ \text{ in }\ \calH\ \ \
  \Rightarrow\ \ \ 
  E_hv_h\to Ev\ \ \text{ in }\ \calE\ \ \text{ with }\
   E\in{\rm Lin}(\calU,\calE)\,.
\end{align}
Here $V_h\subset\calH$ denotes a finite-dimensional space.
Let us note that, in general, we admit a ``non-conformal'' situation that
$V_h\not\subset\calU$ and $E_h\not\in{\rm Lin}(\calU,\calE)$. This allows
for a non-conformal approximation of $v$-variable typically used in
computational implementation of elastodynamics, cf.\ also Section~\ref{sec-numer} below.

The energetics of this scheme can be obtained by imitating
\eqref{energy-test}--\eqref{energy}. More specifically, we test the
particular equations/inclusion in \eqref{suggestion} respectively as
follows: \eqref{suggestion-1} by
$\frac12\varPhi_\varSigma'(\varSigma_{\tau h}^{\KK},\COLOR{z_{\tau h}^{k+1}})
+\frac12\varPhi_\varSigma'(\varSigma_{\tau h}^{\KKK},\COLOR{z_{\tau h}^{k}})$,
then \eqref{suggestion-3} by $\frac{z_{\tau h}^{k+1}-z_{\tau h}^{k}}\tau$, and
eventually the average of \eqref{suggestion-2} at the level $k{+}1$ and ${k}$
by $v_{\tau h}^{k}$. Using that $\varPhi(\cdot,z)$ and $\varPhi(\Sigma,\cdot)$
are quadratic as assumed in \eqref{ansatz}, we have
\begin{subequations}\label{test0}
  \begin{align}\nonumber
&\bigg\langle\frac{\varPhi_\varSigma'(\varSigma_{\tau h}^{\KK},\COLOR{z_{\tau h}^{k+1}})
  +\varPhi_\varSigma'(\varSigma_{\tau h}^{\KKK},\COLOR{z_{\tau h}^{k}})}2,
    \frac{\varSigma_{\tau h}^{\KK}\!-\varSigma_{\tau h}^{\KKK}}\tau\bigg\rangle
    \\&\nonumber\quad=
\bigg\langle\frac{\varPhi_\varSigma'(\varSigma_{\tau h}^{\KK},z_{\tau h}^{k})
  +\varPhi_\varSigma'(\varSigma_{\tau h}^{\KKK},z_{\tau h}^{k})}2,
    \frac{\varSigma_{\tau h}^{\KK}\!-\varSigma_{\tau h}^{\KKK}}\tau\bigg\rangle
    \\&\nonumber\qquad\qquad\qquad\qquad\qquad\qquad\qquad
    \COLOR{+}\frac\tau2\bigg\langle\frac{\varPhi_\varSigma'(\COLOR{\varSigma_{\tau h}^{\KK},z_{\tau h}^{k+1}})
  -\varPhi_\varSigma'(\COLOR{\varSigma_{\tau h}^{\KK},z_{\tau h}^{k}})}\tau,
    \frac{\varSigma_{\tau h}^{\KK}\!-\varSigma_{\tau h}^{\KKK}}\tau\bigg\rangle
    \\&
    \quad=
\frac{\varPhi(\varSigma_{\tau h}^{\KK},z_{\tau h}^{k})
  {-}\varPhi(\varSigma_{\tau h}^{\KKK},z_{\tau h}^{k})}\tau
  \COLOR{+}\frac\tau2\bigg\langle
  \frac{\varPhi_\varSigma'(\COLOR{\varSigma_{\tau h}^{\KK},z_{\tau h}^{k+1}})
  -\varPhi_\varSigma'(\COLOR{\varSigma_{\tau h}^{\KK},z_{\tau h}^{k}})}\tau,
    \frac{\varSigma_{\tau h}^{\KK}\!-\varSigma_{\tau h}^{\KKK}}\tau\bigg\rangle\,,
\label{test01}
\intertext{where we used also \eqref{suggestion-1}, and}
&\bigg\langle\varPhi_z'\Big(\varSigma_{\tau h}^{\KK},
\frac{z_{\tau h}^{k+1}{+}z_{\tau h}^{k}}{2}\Big),
\frac{z_{\tau h}^{k+1}-z_{\tau h}^{k}}\tau\bigg\rangle=
\frac{\varPhi(\varSigma_{\tau h}^{\KK},z_{\tau h}^{k+1})-
\varPhi(\varSigma_{\tau h}^{\KK},z_{\tau h}^{k})}\tau\,.
\end{align}\end{subequations}
Therefore, this test gives
\begin{subequations}\label{test}
  \begin{align}\nonumber
&\frac{\varPhi(\varSigma_{\tau h}^{\KK},z_{\tau h}^{k})
      -\varPhi(\varSigma_{\tau h}^{\KKK},z_{\tau h}^{k})}\tau
    =
\bigg\langle\frac{\varPhi_\varSigma'(\varSigma_{\tau h}^{\KK},\COLOR{z_{\tau h}^{k+1}})
  +\varPhi_\varSigma'(\varSigma_{\tau h}^{\KKK},\COLOR{z_{\tau h}^{k}})}2,\mathfrak CE_h v_{\tau h}^{k}+D_{\tau h}^k\bigg\rangle
\\[-.3em]&\hspace{14em}
\COLOR{-}\frac\tau2\bigg\langle\frac{\varPhi_\varSigma'(\COLOR{\varSigma_{\tau h}^{\KK},z_{\tau h}^{k+1}})
  -\varPhi_\varSigma'(\COLOR{\varSigma_{\tau h}^{\KK},z_{\tau h}^{k}})}\tau,
    \frac{\varSigma_{\tau h}^{\KK}\!-\varSigma_{\tau h}^{\KKK}}\tau\bigg\rangle
\label{test-1}
\,,
\\
&
\varXi
\Big(\frac{z_{\tau h}^{k+1}{-}z_{\tau h}^{k}}\tau\Big)
+\frac{\varPhi(\varSigma_{\tau h}^{\KK},z_{\tau h}^{k+1})-
\varPhi(\varSigma_{\tau h}^{\KK},z_{\tau h}^{k})}\tau\le0\,,
\label{test3}
\\
&\big\langle\TT'\frac{v_{\tau h}^{k+1}{-}v_{\tau h}^{k-1}}{2\tau},v_{\tau h}^{k}
\big\rangle+
\bigg\langle E_h^*\mathfrak C^*
\frac{\varPhi_\varSigma'(\varSigma_{\tau h}^{\KK},\COLOR{z_{\tau h}^{k+1}})
  {+}\varPhi_\varSigma'(\varSigma_{\tau h}^{\KKK},\COLOR{z_{\tau h}^{k}})}2,v_{\tau h}^{k}\bigg\rangle
=\big\langle F_{\tau h}^{k},v_{\tau h}^{k}\big\rangle\,,
\label{test2}
\end{align}\end{subequations}
with $F_{\tau h}^{k}:=\frac12F_{\tau h}^{\KK}{+}\frac12F_{\tau h}^{\KKK}$
and $\varSigma_{\tau h}^{k}:=\frac12\varSigma_{\tau h}^{\KK}{+}\frac12\varSigma_{\tau h}^{\KKK}$.
Let us also note that, if $\varPsi(0)=0$ is assumed, 
the substitution $\widetilde z=0$ into the inequality \eqref{VI}
gives $\varPsi(\frac{z_{\tau h}^{k+1}-z_{\tau h}^{k}}\tau)$
instead of the dissipation rate $\varXi(\frac{z_{\tau h}^{k+1}-z_{\tau h}^{k}}\tau)$
in \eqref{test3},
which is a suboptimal estimate except if $\varPsi$ is degree-1
positively homogeneous.

Summing \eqref{test} up, we enjoy the cancellation of the terms
$\pm\varPhi(\varSigma_{\tau h}^{\KK},z_{\tau h}^{k})$, which is the usual
attribute of the fractional-split scheme. Thus,
using also the simple algebra $\langle\TT'(v_{\tau h}^{k+1}-v_{\tau h}^{k-1}),
v_{\tau h}^{k}\rangle=\langle\TT'v_{\tau h}^{k+1},v_{\tau h}^{k}\rangle
-\langle\TT'v_{\tau h}^{k},v_{\tau h}^{k-1}\rangle$,
we obtain the analog of
\eqref{energy}, namely
\begin{align}\nonumber
&\frac{\langle\TT'v_{\tau h}^{k+1},v_{\tau h}^{k}\rangle
-\langle\TT'v_{\tau h}^{k},v_{\tau h}^{k-1}\rangle}{2\tau}
  +\frac{\varPhi(\varSigma_{\tau h}^{\KK},z_{\tau h}^{k+1})
    -\varPhi(\varSigma_{\tau h}^{\KKK},z_{\tau h}^{k})}\tau
+\varXi\Big(\frac{z_{\tau h}^{k+1}-z_{\tau h}^{k}}\tau\Big)
\\[-.3em]&\nonumber\qquad\qquad\qquad
\le\big\langle F_{\tau h}^{k},v_{\tau h}^{k}\big\rangle
+
\bigg\langle\frac{\varPhi_\varSigma'(\varSigma_{\tau h}^{\KK},\COLOR{z_{\tau h}^{k+1}})
  +\varPhi_\varSigma'(\varSigma_{\tau h}^{\KKK},z_{\tau h}^{k})}2,D_{\tau h}^k\bigg\rangle
\\[-.1em]&\qquad\qquad\qquad\quad\
\COLOR{-}\frac\tau2\bigg\langle\frac{\varPhi_\varSigma'(\COLOR{\varSigma_{\tau h}^{\KK},z_{\tau h}^{k+1}})
  -\varPhi_\varSigma'(\COLOR{\varSigma_{\tau h}^{\KK},z_{\tau h}^{k}})}\tau,
    \frac{\varSigma_{\tau h}^{\KK}\!-\varSigma_{\tau h}^{\KKK}}\tau\bigg\rangle\,.
 \label{energy-disc}\end{align}
If $\varPsi$ is smooth except possibly at zero, there is even equality
in \eqref{energy-disc}.

Considering some approximate values $\{z_{\tau h}^k\}_{k=0,...,K}$ of the variable
$z$ with $K=T/\tau$, we define the piecewise-constant and the piecewise affine 
interpolants respectively by
\begin{subequations}\label{def-of-interpolants}
\begin{align}\label{def-of-interpolants-}
&&&
\overline{z}_{\tau h}(t)= z_{\tau h}^k,\qquad\ \
\underline{\overline z}_{\tau h}(t)=\frac12z_{\tau h}^k+\frac12z_{\tau h}^{k-1},
\qquad\ \text{ and}
&&
\\&&&\label{def-of-interpolants+}
z_{\tau h}(t)=\frac{t-(k{-}1)\tau}\tau z_{\tau h}^k
+\frac{k\tau-t}\tau z_{\tau h}^{k-1}
&&\hspace*{-7em}\text{for }(k{-}1)\tau<t\le k\tau.
\end{align}\end{subequations}
Similar meaning is implied for $\varSigma_{\tau h}$, $v_{\tau h}$,
$\overline\varSigma_{\tau h}$, $\overline v_{\tau h}$, $\overline F_{\tau h}$, etc. 
%
The discrete scheme \eqref{suggestion} can be written in a ``compact'' form as 
\begin{subequations}\label{suggestion+}
\begin{align}\label{suggestion-1+}
  &
  \DT\varSigma_{\tau h}=
  \mathfrak CE_h \overline v_{\tau h}+\DT G_{\tau h}\ \ \text{ and }\ \
  \DT u_{\tau h}=\overline v_{\tau h},
\\\label{suggestion-3+}
&
\partial\varPsi\big(\DT z_{\tau h}\big)+
\varPhi_z'\big(\overline\varSigma_{\tau h},\underline{\overline z}_{\tau h}\big)\ni0\,,
\\\label{suggestion-2+}
  &
  \TT'\DT v_{\tau h}+E_h^\ast\overline S_{\tau h}
  =\overline F_{\tau h}
  \ \ \text{ with }\ \ \overline S_{\tau h}
  =\mathfrak C^*\varPhi_\varSigma'(\overline\varSigma_{\tau h},\COLOR{\overline z_{\tau h}})\,.
\end{align}\end{subequations}

\section{Numerical stability and convergence}
\label{sect-anal}

Because the 
energy \eqref{leap-frog-energy} involves now also the internal variable, 
the CFL condition \label{CFLold}
becomes
\begin{align}
\exists\,\eta>0\ \forall\,\varSigma_h,z_h,\wt z_h:\ \ 
  \varPhi(\varSigma_h,z_h)\ge\frac{\tau^2}{4{-}\eta}
      \big\langle E^*_h S_h,(\TT')^{-1}  E^*_h S_h\big\rangle_{{\calH}^*\times{\calH}}^{}
      \ \ \text{ with }\ S_h=\mathfrak C^* \varPhi'_{\Sigma} (\varSigma_h,\wt z_h),
  \label{CFLnew}
\end{align}
where
$\varSigma_h$, $z_h$,  and $\wt z_h$ is considered
from the corresponding finite-dimensional subspaces. 
Let us still introduce the Banach space
${\mathcal X}:=\{X\in{\calS}^*;\ E^*{\mathfrak C}^*X\in{\calH}^*\}$.
We further assume
$\mathfrak C\in{\rm Lin}({\calE},{\calS})$
invertible.

\begin{proposition}[Numerical stability.]\label{prop1}
  Let $F$ be constant in time, valued in ${\calH}^*$,
  $G\!\in\!W^{1,1}(I;{\calS})$,
  $u_0\in\mathcal\calU$ so that $\varSigma_0={\mathfrak C}Eu_0\in
  {\calS}$, $v_0\in\calH$, $z_0\in\calZ$, the functionals
  $\mathcal T$, $\varPhi$, and $\varPsi$ be coercive and 
  $\varPhi_\varSigma'(\varSigma,\cdot)$ be Lipschitz continuous uniformly for
  $\varSigma\in\calS$ in the sense
\begin{subequations}\label{coerc+Lipschitz}\begin{align}\nonumber
    &&&\exists\epsilon>0\ p\ge2\ \forall (\varSigma,v,z)\in
    {\calS}{\times}{\calH}{\times}\calZ:
    \\&&&&&\hspace*{-9em}\mathcal T(v)\ge\epsilon\|v\|_{\calH}^2,\ \
    \varPhi(\varSigma,z)\ge\epsilon\|\varSigma\|_{\calS}^2+\epsilon\|z\|_\calZ^2,
    \ \
  \varPsi(z)\ge\epsilon\|z\|_{\calZ_1}^p,&&&&
  \\&&&\label{ass-Phi_S}
  \exists C\ \forall \varSigma\in\calS,\ z\in \calZ:&&
  \big\|\varPhi_\varSigma'(\varSigma,z)\big\|_{{\calS}^*}^{}\le
  C\big(1+\|\varSigma\|_{\calS}^{}+\|z\|_\calZ^{}\big),
  \\&&&\exists\ell\in\R\ \forall \varSigma\in\calS,\ z,\wt z\in \calZ:&&
    \big\|\varPhi_\varSigma'(\varSigma,z)-\varPhi_\varSigma'(\varSigma,\wt z)\big\|_{{\calS}^*}^{}
    \le\ell\|z-\wt z\big\|_{\calZ_1}^{}.
    \label{Lipschitz}
  \end{align}\end{subequations}
  Let also the CFL condition \eqref{CFLnew} hold with
  $\tau>0$ sufficiently small (in order to make the discrete Gronwall
  inequality effective).
Then the following a-priori estimates hold:
\begin{subequations}\label{est}
  \begin{align}\label{est0}
    &\|u_{\tau h}\|_{W^{1,\infty}(I;\mathcal H)}^{}\le C\,,
    \\\label{est1}
    &\|\varSigma_{\tau h}\|_{L^\infty(I;\calS)}^{}\le C\ \ \text{ and }\ \
    \|\DT\varSigma_{\tau h}\|_{L^1(I;\mathcal X^*)}^{}\le C,
    \\&\label{est2}\|v_{\tau h}\|_{L^\infty(I;\calH)}^{}\le C\ \ \text{ and }\ \
    \|\mathcal T'\DT v_{\tau h}\|_{L^\infty(I;\calU^*)}^{}\le C,
    \\&\label{est3}\|z_{\tau h}\|_{L^\infty(I;\calZ)}^{}\le C
    \ \ \text{ and }\ \
    \|\DT z_{\tau h}\|_{L^p(I;\calZ_1)}^{}\le C,
\end{align}\end{subequations}    
\end{proposition}

\noindent{\it Proof.}
The energy 
{imbalance} that we have here is \eqref{energy-disc} which can be re-written as
 \begin{align}\nonumber
&\frac{\mathfrak E_h^{\KK}{-}\mathfrak E_h^{\KKK}}{\tau}+\varXi
\Big(\frac{z_{\tau h}^{k+1}{-}z_{\tau h}^{k}}\tau\Big)
\le\big\langle F_{\tau h}^{k},v_{\tau h}^{k}\big\rangle_{{\calH}^*\times{\calH}}^{}
+\bigg\langle\frac{\varPhi_\varSigma'(\varSigma_{\tau h}^{\KK},\COLOR{z_{\tau h}^{k+1}})
  +\varPhi_\varSigma'(\varSigma_{\tau h}^{\KKK},z_{\tau h}^{k})}2,D_{\tau h}^k\bigg\rangle_{{\calS}^*\times{\calS}}^{}
\\[-.1em]&\qquad\qquad\qquad\qquad\qquad\qquad
\COLOR{-}\frac\tau2\bigg\langle\frac{\varPhi_\varSigma'(\COLOR{\varSigma_{\tau h}^{\KK},z_{\tau h}^{k+1}})
  -\varPhi_\varSigma'(\COLOR{\varSigma_{\tau h}^{\KK},z_{\tau h}^{k}})}\tau,
    \frac{\varSigma_{\tau h}^{\KK}\!-\varSigma_{\tau h}^{\KKK}}\tau
   \bigg\rangle_{{\calS}^*\times{\calS}}^{}
 \label{en1}\end{align}
 with an analog of the energy \eqref{leap-frog-energy}, namely
 \begin{align}
   \mathfrak E_h^{\KK}  =
   \frac12\langle\TT'v_{\tau h}^{k+1},v_{\tau h}^{k}\rangle_{{\calH}^*\times{\calH}}^{}
 +\varPhi(\varSigma_{\tau h}^{\KK},z_{\tau h}^{k+1}).
  \label{en2}
\end{align}
 We need to show that $\mathfrak E_h^{k+1} $ is indeed 
 {a sum of the kinetic and the stored energies}
 { at least up to some positive coefficients}. To do so, like e.g.\
\cite[Lemma~4.2]{Scar04ETNP} or \cite[Sect.~6.1.6]{Tsog99}, let us write
\begin{align}
  \langle\TT'v_{\tau h}^{k+1},v_{\tau h}^{k}\rangle 
  &=  \Big\langle\TT' \frac{v_{\tau h}^{k+1}+v_{\tau h}^{k} }{2} ,\frac{v_{\tau h}^{k+1}+v_{\tau h}^{k} }{2}\Big\rangle
  -  \Big\langle\TT' \frac{v_{\tau h}^{k+1}-v_{\tau h}^{k} }{2} ,\frac{v_{\tau h}^{k+1}-v_{\tau h}^{k} }{2}\Big\rangle
  \nonumber \\\nonumber
  &=  \Big\langle \TT' \frac{v_{\tau h}^{k+1}+v_{\tau h}^{k} }{2} ,\frac{v_{\tau h}^{k+1}+v_{\tau h}^{k} }{2}\Big\rangle
  \\
  &\ - \frac{\tau^2}{4}\big\langle  E^*_h\big(\mathfrak C^*\varPhi'_\Sigma(\varSigma^{\KK}_{\tau h},\COLOR{z_{\tau h}^{k+1}}){-}F_{\tau h}^{\KK}\big),
  (\TT')^{-1} E^*_h\big(\mathfrak C^* \varPhi'_\Sigma( \varSigma^{\KK}_{\tau h},\COLOR{z_{\tau h}^{k+1}}){-}F_{\tau h}^{\KK}\big)\big\rangle\,,
\label{en3}
\end{align}
where all the duality pairings are between ${\calH}^*$ and
${\calH}$; 
here also \eqref{suggestion-2} has been used.
Thus, using also $\TT(v)=\frac12\langle\TT'v,v\rangle$,
we can write the energy \eqref{en2} as
\begin{align}\nonumber
  \mathfrak E_h^{\KK}&=
\TT(v_{\tau h}^{k+1/2})
+a_{\tau h}^{\KK}\varPhi(\varSigma_{\tau h}^{\KK},z_{\tau h}^{k+1})
+\frac{\tau^2}{2}\big\langle(\TT')^{-1}E^*_h\mathfrak C^*
\varPhi'_\Sigma(\varSigma^{\KK}_{\tau h},\COLOR{z_{\tau h}^{k+1}}),F_{\tau h}^{\KK}\big\rangle
- \frac{\tau^2}{4}\|F_{\tau h}^{\KK}\|_{\calH}^2   \\
  & \text{ with }\ \ a_{\tau h}^{\KK}:=
  1-\frac{\tau^2}{4}  \frac{  \big\langle  E^*_h \mathfrak C^* \varPhi'_\Sigma( \varSigma^{\KK}_{\tau h},\COLOR{z_{\tau h}^{k+1}}), (\TT')^{-1} E^*_h \mathfrak C^* \varPhi'_\Sigma( \varSigma^{\KK}_{\tau h},\COLOR{z_{\tau h}^{k+1}})\big\rangle}{\varPhi(\varSigma_{\tau h}^{\KK},z_{\tau h}^{k+1}) }\ge\eta
  \label{en4}\end{align}
and with $v_{\tau h}^{k+1/2}:=\frac12v_{\tau h}^{k+1}+\frac12v_{\tau h}^{k}$.
The energy $\mathfrak E_h^{\KK}$
        {yields a-priori estimates} if the
coefficient $a_{\tau h}^k$ 
is non negative, which is just ensured by our CFL
  condition
%
  \eqref{CFLnew} used for $\varSigma_h=\varSigma_{\tau h}^{\KK}$, $z_h=z_{\tau h}^{k+1}$and $\wt z_h=z_{\tau h}^k$. 
  Here $\eta>0$ is just from \eqref{CFLnew}. 


Altogether, summing \eqref{en1} for $k=0,...,l-1\le T/\tau-1$ and using \eqref{en4},
we obtain the estimate 
 \begin{align}\nonumber
   &
   \epsilon\bigg(
   \big\|v_{\tau h}^{l-1/2}\big\|_{\calH}^2
   +a_{\tau h}^{l-1}\big\|\varSigma_{\tau h}^{l}\big\|_{\calS}^2
   +a_{\tau h}^{l-1}\big\|z_{\tau h}^{l}\big\|_\calZ^2
   +\tau\sum_{k=0}^{l-1}
   \Big\|\frac{z_{\tau h}^{k+1}{-}z_{\tau h}^{k}}\tau\Big\|_{\calZ_1}^p\bigg)
\\&\nonumber\quad
\le
\frac{\tau^2}{4}\|F_{\tau h}^{l}\|_{\calH}^2
-\frac{\tau^2}{2}\big\langle(\TT')^{-1}E^*_h\mathfrak C^*
\varPhi'_\Sigma(\varSigma^{l}_{\tau h},\COLOR{z_{\tau h}^{l}}),F_{\tau h}^{l}\big\rangle
-\frac{\tau^2}{2}\big\langle(\TT')^{-1}E^*_h\mathfrak C^*
\varPhi'_\Sigma(\varSigma^{0}_{\tau h},\COLOR{z_{\tau h}^{0}}),F_{\tau h}^{0}\big\rangle
\\[-.5em]&\nonumber\quad\
   +\TT(v_{\tau h}^{-1/2})
   +a_{\tau h}^{0}\varPhi(\varSigma_{\tau h}^{0},z_{\tau h}^{0})
+\tau\sum_{k=0}^{l-1}\bigg(\big\langle F_{\tau h}^{k},v_{\tau h}^{k}\big\rangle
\\[-.5em]&
\quad\
   +\frac12\big\|\varPhi_\varSigma'(\varSigma_{\tau h}^{\KK},\COLOR{z_{\tau h}^{k+1}})
  +\varPhi_\varSigma'(\varSigma_{\tau h}^{\KKK},z_{\tau h}^{k})\big\|_{{\calS}^*}
\big\|D_{\tau h}^k\big\|_{{\calS}}^{}
+\frac\tau2\ell\Big\|\frac{z_{\tau h}^{k}{-}z_{\tau h}^{k-1}}\tau\Big\|_{\calZ_1}
\Big\|
\frac{\varSigma_{\tau h}^{\KK}\!{-}\varSigma_{\tau h}^{\KKK}}\tau\Big\|_{{\calS}}^{}\bigg)\,,
 \label{overall}\end{align}
 where $\epsilon$, $p$, $\ell$ and $a_{\tau h}^{l}$ come from
 \eqref{coerc+Lipschitz} and \eqref{en4}.
Using \eqref{ass-Phi_S}, we estimate 
$\|\varPhi_\varSigma'(\varSigma_{\tau h}^{\KK},z_{\tau h}^{k})
+\varPhi_\varSigma'(\varSigma_{\tau h}^{\KKK},z_{\tau h}^{k})\|_{{\calS}^*}
\le C(1+\|\varSigma_{\tau h}^{\KK}\|_{\calS}^2
+\|z_{\tau h}^{\KK}\|_\calZ^2)$ and then use
the summability of $\|D_{\tau h}^k\big\|_{{\calS}}^{}$ needed
for the discrete Gronwall inequality; here the assumption
$\DT G\in L^1(I;{\calS})$ is needed.
 The last term in \eqref{overall} is
 to be estimated by the H\"older inequality as
\begin{align}
\frac\tau2\ell\Big\|\frac{z_{\tau h}^{k}{-}z_{\tau h}^{k-1}}\tau\Big\|_{\calZ_1}
\Big\|
\frac{\varSigma_{\tau h}^{\KK}\!{-}\varSigma_{\tau h}^{\KKK}}\tau\Big\|_{\calS}^{}
\le\frac\epsilon2\Big\|\frac{z_{\tau h}^{k}{-}z_{\tau h}^{k-1}}\tau\Big\|_{\calZ_1}^p
+C_{p,\epsilon,\ell}^{}\big(1+
\|\varSigma_{\tau h}^{\KK}\|_{\calS}^2
+\|\varSigma_{\tau h}^{\KKK}\|_{\calS}^2\big)\end{align}
with some $C_{p,\epsilon,\ell}^{}$ depending on $p$, $\epsilon$, and $\ell$. Here
we needed $p\ge2$; note that this is related with the specific explicit time
discretisation due to the last term in \eqref{energy-disc} but not with the
problem itself. Then we use
the discrete Gronwall inequality
to obtain the former estimates in (\ref{est}b,c) and the estimates
(\ref{est}a,d).
The usage of the mentioned discrete Gronwall inequality is however a bit
tricky because of the the term $\|v_{\tau h}^{l-1/2}\big\|_{\calH}^2$
on the left-hand side of \eqref{overall} while there is
$v_{\tau h}^{k}$ instead of $v_{\tau h}^{k-1/2}$. To cope with it,
we rely on $F$ constant (as assumed) and,
proving the estimate
for $l=1$, we sum up \eqref{overall} for $l+1$ and $l$ to see
$\big\langle F_{\tau h}^{k},v_{\tau h}^{k-1/2}\big\rangle$ also
on the right-hand side. 
\COMMENT{TO CHECK DETAILS!!}

The equation $\DT\varSigma_{\tau h}=\mathfrak CE_h\overline v_{\tau h}
+\DT G_{\tau h}$ gives the latter estimate in
\eqref{est1} by estimating
\begin{align}\nonumber
  \int_0^T\!\!
  \big\langle\DT\varSigma_{\tau h},X\big\rangle_{\mathcal X^*\times\mathcal X}\,\d t
&=\int_0^T\!\!\big\langle\mathfrak CE_h\overline v_{\tau h}
+\DT G_{\tau h},X\big\rangle_{\mathcal X^*\times\mathcal X}\,\d t
\\&=\int_0^T\!\!\big\langle\overline v_{\tau h},E_h^*\mathfrak C^*X\big\rangle_{\calH\times\calH^*}\,\d t
+\int_0^T\!\!\big\langle\DT G_{\tau h},X\big\rangle_{\mathcal X^*\times\mathcal X}\,\d t
\label{dual-est1}\end{align}
for $X\in L^\infty(I;\mathcal X)$ and using also the already proved boundedness
of $\overline v_{\tau h}$ in $L^\infty(I;\calH)$ and the assumed
boundedness of $E_h$ uniform in $h>0$;
here we used also that $\DT\varSigma_{\tau h}(t)\in\calS\subset\mathcal X^*$.

Eventually, the already obtained estimates \eqref{ass-Phi_S} 
give $\varPhi_\varSigma'(\overline\varSigma_{\tau h},\COLOR{\overline z_{\tau h}})$
bounded in $L^\infty(I;\calS^*)$. Therefore 
$\overline S_{\tau h}=\mathfrak C^*\varPhi_\varSigma'(\overline\varSigma_{\tau h},\COLOR{\overline z_{\tau h}})$ is bounded in $L^\infty(I;\calE^*)$, hence
$E_h^*\overline S_{\tau h}$ is bounded in $L^\infty(I;\calU^*)$,
so that $\mathcal T'\DT v_{\tau h}=\overline F_{\tau h}-E_h^*\overline S_{\tau h}$
gives the latter estimate in \eqref{est2}.
$\hfill\Box$


\begin{proposition}[Convergence.]\label{prop2}
  Let \eqref{z-ansatz} and \eqref{Eh->E} hold, all the involved Banach spaces
  be separable, 
  and the assumptions of Proposition~\ref{prop1} hold. Moreover, let
\begin{subequations}\label{ass-Phi'}\begin{align}\nonumber
    &\forall z\In\calZ:\ \varPhi_\varSigma'(\cdot,z)\ \text{continuous linear, and }
    \varPhi_\varSigma'
    :\calS{\times}\calZ_0\to{\rm Lin}(\calS,\calS^*)\ \text{continuous}
  \\&\label{ass-Phi-Sigma}\qquad\qquad\qquad\qquad\qquad
\text{or }\ \varPhi_\varSigma':\calS\times\calZ\to\calS^*\
  \text{ is continuous linear,}
  \\\nonumber
  &\forall z\In\calZ:\ [\varPhi_0]_z'(\cdot,z)\ \text{continuous linear, and }
      [\varPhi_0]_z'
      :\calS{\times}\calZ_0\to{\rm Lin}(\calZ_0,\calZ_1)\ \text{continuous}
  \\&\qquad\qquad\qquad\qquad\qquad\text{or }\ [\varPhi_0]_z':\calS\times\calZ\to\calZ_1^*\
  \text{  is continuous linear, and}
  \label{ass-Phi0-z}
  \\&\varPhi_1':\calZ\to\calZ^*\ \text{ is linear continuous},
  \end{align}\end{subequations}
for some Banach space $\calZ_0$ into which $\calZ$ is embedded compactly,
where $\varPhi_0$ and $\varPhi_1$ are from \eqref{z-ansatz}.
Then there is a selected subsequences, again denoting
  $\{(u_{\tau h},\varSigma_{\tau h},v_{\tau h},z_{\tau h})\}_{\tau>0}$
  converging weakly* in the topologies indicated in the estimates \eqref{est}
  to some $(u,\varSigma,v,z)$.
  Moreover, any $(u,\varSigma,v,z)$ obtained as such a limit is a weak
  solution according Definition~\ref{def}.
  \end{proposition}

\noindent{\it Proof.}
By the Banach selection principle, we can select the weakly* converging
subsequence as claimed; here the separability of the involved Banach spaces 
is used.

Referring to the compact embedding $\calZ\subset\calZ_0$ used in
the former option in (\ref{ass-Phi'}a,b) and relying on a generalization
the Aubin-Lions compact-embedding theorem with
$\COLOR{\DT{\overline z}_{\tau h}}$
being bounded in the space of the $\calZ_1$-valued measures on $I$,
cf.\ \cite[Corollary 7.9]{Roub13NPDE}, we have $\COLOR{\overline z_{\tau h}}\to z$
strongly in $L^{r}(I;\calZ_1)$ for any $1\le r<+\infty$. 

Further, we realize that the approximate solution satisfy identities/inequality
analogous to what is used in Definition~\ref{def}.
In view of \eqref{weak-v}, the equations \eqref{suggestion-2+} now
  means
  \begin{subequations}\label{weak}\begin{align}
&\label{weak-v-disc}
\int_0^T\!\!\big\langle \varPhi_\varSigma'(\overline\varSigma_{\tau h},\COLOR{\overline z_{\tau h}}),\mathfrak C E_h\wt v\big\rangle_{{\calS}^*\times{\calS}}^{}
    -\big\langle\mathcal T'v_{\tau h},\DT{\wt v}\big\rangle_{{\calH}^*\times{\calH}}^{}
    \,\d t=\big\langle\mathcal T'v_0,\wt v(0)\big\rangle_{{\calH}^*\times{\calH}}^{}
    +\int_0^T\!\!\big\langle F_{h},\wt v\big\rangle_{{\calH}^*\times{\calH}}^{}\,\d t
    \intertext{for any $\wt v\in C^1(I;{\calH})
      $ valued in $V_h$ \COMMENT{OK??} and 
      with $\wt v(T)=0$.
 In view of \eqref{weak-z}, the inclusion \eqref{suggestion-3+} means}
 &\nonumber
   \int_0^T\!\!\varPsi(\wt z)+\big\langle[\varPhi_0]_z'(\overline\varSigma_{\tau h},\underline{\overline z}_{\tau h}),\wt z{-}
    \DT z_{\tau h}\big\rangle_{\calZ_1^*\times\calZ_1}^{}\!
    +\big\langle\varPhi_1'\underline{\overline z}_{\tau h},\wt z\big\rangle_{\calZ^*\times\calZ}^{}
    \,\d t+\varPhi_1(z_0)
    \\[-.4em]&\hspace*{23em}
    \ge\varPhi_1(z_{\tau h}(T))+
    \!\int_0^T\!\!\varPsi(\DT z_{\tau h})\,\d t\,.
    \label{weak-z-disc}
\end{align}\end{subequations}
This is completed by \eqref{suggestion-1+}. 

It is further important that the equations in
\eqref{suggestion-1+} and the first equation in \eqref{suggestion-2+}
are linear, so that the weak convergence is sufficient for the
limit passage there. In particular, we use \eqref{Eh->E} and the
Lebesgue dominated-convergence theorem.

As to the weak convergence of \eqref{suggestion-1+} integrated in time
towards \eqref{suggestion-1} integrated in time, i.e.\ towards $\varSigma
=\mathfrak CEu+G$ as used in Definition~\ref{def}, we need to prove that
\begin{align}
\int_0^T\!\!\big\langle\varSigma_{\tau h}-G_{\tau h},X\big\rangle_{\calS\times\calS^*}
-\big\langle u_{\tau h},E_h^*\mathfrak C^*X\big\rangle_{\calH\times\calH^*}\d t
\to\int_0^T\!\!\big\langle\varSigma{-}G,X\big\rangle_{\calS\times\calS^*}
-\big\langle u,E^*\mathfrak C^*X\big\rangle_{\calH\times\calH^*}\d t
\label{conv-S=CEu}\end{align}
for any $X\in\calS^*$.
By \eqref{Eh->E}, we have also $E_h^*S\to E^*S$ in $\calH$ for any
$S\in\calE^*$, in particular for $S=\mathfrak C^*X(t)$.
Thus certainly $E_h^*\mathfrak C^*X\to E^*\mathfrak C^*X$ in $L^1(I;\calH)$
strongly. Using the weak* convergence $u_{\tau h}\to u$ in
$L^\infty(I;\calH)$, we obtain \eqref{conv-S=CEu}.
Moreover, in the limit
$Eu=\mathfrak C^{-1}(\varSigma-G)\in L^\infty(I;\calE)$ so that
$u\in L^\infty(I;\calU)$.

For the limit passage in \eqref{weak-v-disc}, we also use
$\varPhi_\varSigma'(\overline\varSigma_{\tau h},\COLOR{\overline z_{\tau h}})\to
\varPhi_\varSigma'(\varSigma,z)$ weakly* in $L^\infty(I;\calS^*)$
because $\varPhi_\varSigma'$ is continuous in the
(weak$\times$strong,weak)-mode, cf.~\eqref{ass-Phi-Sigma},
and because of the mentioned strong convergence of
$\overline z_{\tau h}\to z$.


Furthermore, we need to show the convergence 
$[\varPhi_0]_z'(\overline\varSigma_{\tau h},\underline{\overline z}_{\tau h})\to
[\varPhi_0]_z'(\varSigma,z)$.
For this, we use again the mentioned generalized Aubin-Lions theorem
to have the strong convergence $\underline{\overline z}_{\tau h}\to z$
in $L^{r}(I;\calZ_1)$ for any $1\le r<+\infty$ and then
the continuity of $[\varPhi_0]_z'$ in the 
(weak$\times$strong,weak)-mode, cf.~the former option in \eqref{ass-Phi0-z}.
The limit passage of \eqref{weak-z-disc} towards
\eqref{weak-z} then uses also the weak lower semicontinuity of $\varPhi_1$ and
the weak convergence $z_{\tau h}(T)\to z(T)$ in $\calZ$; here for this
pointwise convergence in all time instants $t$ and in particular in $t=T$,
we also used that we have some information about $\DT z_{\tau h}$,
cf.\ \eqref{est3}.

So far, we have relied on the former options in (\ref{ass-Phi'}a,b) and the
Aubin-Lions compactness argument as far as the $z$-variable concerns.
If $\varPhi$ is quadratic (as e.g.\ in the examples in
Sects.\,\ref{sect-plast}--\ref{sect-poro} below), we can use the latter
options in (\ref{ass-Phi'}a,b) and simplify the above arguments, relying
merely on the weak convergence $\overline z_{\tau h}\to z$ and
$\overline{\underline z}_{\tau h}\to z$.
$\hfill\Box$

\begin{remark}[Alternative weak formulation]
  \upshape
  Here, we used the weak formulation of \eqref{IVP2z+} containing the term
$\langle\varPhi_z'(\varSigma,z),\DT z\big\rangle$
which often does not have a good meaning since $\DT z$
may not be enough regular in some applications. This term is thus
eliminated
by substituting it, after integration over the time interval, by
$\varPhi(\varSigma(T),z(T))
-\int_0^T\langle\varPhi_\varSigma'(\varSigma,z),\DT\varSigma\big\rangle\,\d t
-\varPhi(\varSigma_0,z_0)$ or even rather by
$\varPhi(\varSigma(T),z(T))
-\int_0^T\langle\varPhi_\varSigma'(\varSigma,z),\mathfrak CEv\big\rangle\,\d t
-\varPhi(\varSigma_0,z_0)$. Here, however, it would bring even more difficulties
because we would need to prove
a strong convergence of
$\varPhi_\varSigma'(\varSigma,z)$, or of $\DT\varSigma$, or $\mathfrak CEv$ in
our explicit-discretisation scheme, which seems not easy. 
\end{remark}

\begin{remark}[Nonquadratic $\varPhi(\varSigma,\cdot)$.]\label{rem-nonquadratic}
\upshape
Some applications use such $\varPhi(\varSigma,\cdot)$ which is not quadratic.
This is still consistent with the explicit leap-frog-type discretisation if,
instead of $\varPhi_z'(\varSigma,z)$,
we consider an abstract difference quotient
$\varPhi_z^\circ(\varSigma,z,\widetilde z)$ with the properties
\begin{align}\label{quotient}
  \varPhi_z^\circ(\varSigma,z,z)=\varPhi_z'(\varSigma,z)
  \ \ \text{ and }\ \ \big\langle\varPhi_z^\circ(\varSigma,z,\widetilde z),
  z{-}\widetilde z\big\rangle
  =\varPhi(\varSigma,z)-\varPhi(\varSigma,\widetilde z)\,,
\end{align}
cf.\ \cite{RouPan17ECTD}. Then, instead of
$\varPhi_z'(\varSigma_{\tau h}^{\KK},\frac{z_{\tau h}^{k+1}{+}z_{\tau h}^{k}}{2})$ in
\eqref{suggestion-3}, to write
$\varPhi_z^\circ(\varSigma_{\tau h}^{\KK},z_{\tau h}^{k+1},z_{\tau h}^{k})$.
\end{remark}

\begin{remark}[State-dependent dissipation.]\label{rem-dissip}
\upshape
The generalization of $\varPsi$ dependent also on $z$ or even on $(\varSigma,z)$
is easy. Then $\partial\varPsi$ is to be replaced by the partial subdifferential
$\partial_{\DT z}\varPsi$ and \eqref{suggestion-3}
should use $\varPsi(\varSigma_{\tau h}^{\KK},z_{\tau h}^{k},\cdot)$
instead of $\varPsi(\cdot)$.
\end{remark}

\begin{remark}[Spatial numerical approximation]\upshape
  From the coercivity of the stored energy $\varPhi$, we have
  $\varSigma_{\tau h}^k\in\calS$ for any $k=0,1,...$ and thus,
  from \eqref{suggestion-1},
  $E_hv_{\tau h}^k\in\calE$ so that $v_{\tau h}^k\in\calU$,
  although the limit $v$ cannot be assumed valued in $\calU$
  in general. Similarly, from \eqref{suggestion-2}, one can read that
  $E_h^*S_{\tau h}^k\in\calH$ although this cannot be expected in
  the limit in general. Anyhow, on the time-discrete level, one can use the
  FEM discretisation similarly as in the linear elastodynamics where
  regularity can be employed, cf.\ \cite{BeJoTs02NFMF,BeJoTs01MFES,Tsog99} for a mixed finite\COLOR{-element} method and \cite{ChLaQi15SDGM} for the more recently developed staggered discontinuous Galerkin method for elastodynamics. 
  \COMMENT{I added $\mbox{\cite{ChLaQi15SDGM}}$ if you know of other relevant references please add??? }
\end{remark}

  \section{Particular examples}\label{sec-exa}

We present four examples from continuum mechanics of deformable
bodies at small strains of different characters to illustrate applicability
of the ansatz \eqref{ansatz} and the above discretisation scheme.
Various combinations of these examples are possible, too, covering
thus a relatively wide variety of models.

We use a standard notation concerning function spaces.
Beside the Lebesgue $L^p$-spaces, we denote by $H^k(\varOmega;\R^n)$
the Sobolev space of functions whose distributional derivatives
are from $L^2(\varOmega;\R^{n\times d^k})$.

\subsection{Plasticity or creep}\label{sect-plast}

The simplest example with quadratic stored energy and local dissipation potential
is the model of plasticity or creep. The internal variable is then the plastic
strain $\pi$, valued in the set of symmetric trace-free matrices
$\R_{\rm dev}^{d\times d}=\{P\in\R^{d\times d};\ P^\top=P,\ {\rm tr}\,P=0\}$.
For simplicity, we consider only homogeneous Neumann or Dirichlet boundary
conditions, so that simply $E=e(u)$ and $\mathfrak C=\bbC$. 
The 
stored energy in terms of strain $e(u)$ is
\begin{align}
\mathscr{W}(u,\pi)=\int_\varOmega\frac12\mathbb C(e(u){-}\pi){:}
(e(u){-}\pi)\,\d x
\,,\label{W-creep}
\end{align}
which is actually a function of the elastic strain
$e_{\rm el}=e{-}\pi$. The additive decomposition $e(u)=e_{\rm el}{+}\pi$ is
referred to as Green-Naghdi's \cite{GreNag65GTEP} decomposition.
This energy leads to
\begin{align}
  \varPhi(\sigma,\pi)=\int_\varOmega\frac12\mathbb C^{-1}\sigma{:}\sigma
  -\sigma{:}\pi+\frac12\mathbb C\pi{:}\pi\,\d x
  \ \ \ \text{ with }\ \sigma=\mathbb C e(u)\,.
  \label{Phi-creep}\end{align} 
Let us note that $\varPhi_\sigma'(\sigma,\pi)=
\mathbb C^{-1}\sigma-\pi=e{-}\pi$, i.e.\ the elastic strain $e_{\rm el}$, and that
the proto-stress $\varSigma=\sigma$ is indeed different from the actual stress
$\sigma-\bbC\pi$.

The dissipation potential is standardly chosen as
\begin{align}
  \varPsi(\DT\pi)=\int_\varOmega\sigma_{_{\rm Y}}|\DT\pi|
  +\frac12\bbD\DT\pi{:}\DT\pi\,\d x
\end{align}
with $\sigma_{_{\rm Y}}\ge0$ a prescribed yield stress and $\bbD$ a positive
semidefinite viscosity tensor. The dissipation rate is then
$\varXi(\DT\pi)=\int_\varOmega\sigma_{_{\rm Y}}|\DT\pi|
  +\bbD\DT\pi{:}\DT\pi\,\d x$.
For $\bbD>0$ and $\sigma_{_{\rm Y}}=0$, we obtain mere creep model or, in
other words, the linear \emph{viscoelastic} model in the
\emph{Maxwell rheology}.
For both $\bbD>0$ and $\sigma_{_{\rm Y}}>0$, we obtain \emph{viscoplasticity}.
For $\bbD=0$ and $\sigma_{_{\rm Y}}>0$,
we would obtain the rate-independent (perfect) plasticity but
our Proposition~\ref{prop1} does not cover this case (i.e.\ $p=1$ is not
admitted).

The functional setting is $\calH=L^2(\varOmega;\R^d)$,
$\calE=\calS=\calZ=\calZ_1
=L^2(\varOmega;\R_{\rm sym}^{d\times d})$ where
$\R_{\rm sym}^{d\times d}$ denotes symmetric $(d{\times}d)$-matrices.
Thus $\calU:=\{v\In L^2(\varOmega;\R^d);\ e(v)\In L^2(\varOmega;\R^{d\times d})\}
=H^1(\varOmega;\R^d)$ by Korn's inequality.

A modification of the stored energy  models 
an {\it isotropic hardening}, enhancing \eqref{W-creep} as
\begin{align}
\mathscr{W}(u,\pi)=\int_\varOmega\frac12\mathbb C_1(e(u){-}\pi){:}
(e(u){-}\pi)+\frac12\mathbb C_2\pi{:}\pi\,\d x
\label{W-creep+}
\end{align}
so that the energy $\varPhi$ from \eqref{Phi-creep} is modified as
\begin{align}
\varPhi(\sigma,\pi)=\int_\varOmega\frac12\mathbb C_1^{-1}\sigma{:}\sigma
-\sigma{:}\pi+\frac12(\mathbb C_1{+}\mathbb C_2)\pi{:}\pi\,\d x\,.
\label{Phi-creep+}
\end{align}
In the pure creep variant $\sigma_{_{\rm Y}}=0$, this is actually the
{\it standard linear solid}
(in a so-called Zener form), considered together with the leap-frog
time discretisation in \cite{BeEzJo04MFEA}. The isochoric constraint
${\rm tr}\,\pi=0$ can then be avoided, assuming that $\mathbb C_2$ is
positive definite.

All these models
lead to a flow rule which is
localized on each element when 
 an element-wise constant approximation of $\pi$ is used, 
and the combination with the explicit discretisation
of the other equations leads to a very fast computational procedure.

Another modification for {\it gradient plasticity} by adding terms
$\frac12\kappa|\nabla\pi|^2$ into the stored energy is easily
possible, too. This modification uses
$\calZ=H^1(\varOmega;\R_{\rm sym}^{d\times d})$ and
\eqref{z-ansatz} with
$\varPhi_1(z)=\int_\varOmega\frac12\kappa|\nabla\pi|^2$ and
makes, however, the flow rule
nonlocal but at least on can benefit from that the usual space discretisation
of the proto-stress $\sigma$ uses the continuous piecewise smooth elements
which allows for handling gradients $\nabla\pi$ if used consistently also for
$\pi$. 

For the quasistatic variant of this model, we refer to the classical
monographs \cite{HanRed99PMTN,Tema85MPP}, while the dynamical model with
$\bbD=0$ is e.g.\ in \cite[Sect.5.2]{MieRou15RIST}.

Noteworthy, all these models bear time regularity if the loading is smooth
and initial conditions regular enough, which can be advantageously reflected
in space FEM approximation, too.

\subsection{Poroelasticity in isotropic materials}\label{sect-poro}

 Another example with quadratic stored energy but less trivial dissipation potential is a
 saturated Darcy or Fick flow of a diffusant in porous media, e.g.\ water in porous elastic
 rock or concrete, or a solvent in elastic polymers. The most simple model is the
 classical Biot model \cite{Biot41GTTS}, capturing effects as swelling or seepage.
In one-component flow, 
 the internal variable is then the scalar-valued diffusant content (or concentration) denoted
 by $\zeta$.

 As in the previous section~\ref{sect-plast}, we consider only Neumann or
 Dirichlet boundary conditions,
 so that $E=e(u)$. Here we use the orthogonal decomposition $e={\rm sph}\,e+{\rm dev}\,e$
 with the spherical (volumetric) part ${\rm sph}\,e:=({\rm tr}\,e)\mathbb I/d$
 and the deviatoric part ${\rm dev}\,e$ and 
 confine ourselves to isotropic materials where the elastic-moduli tensor
 $\mathbb C_{ijkl}=K\delta_{ij}\delta_{kl}+G(\delta_{ik}\delta_{jl}
 +\delta_{il}\delta_{jk}-2\delta_{ij}\delta_{kl}/d)$ with
 $K$ the bulk modulus and $G$ the shear modulus (=\,the second Lam\'e
 constant), which is the standard notation hopefully without any
   confusion with the notation used in \eqref{IVP+}.
 Then the proto-stress $\varSigma=\sigma=\mathbb Ce=K{\rm sph}\,e+2G{\rm dev}\,e$.
 In particular, ${\rm sph}\,\sigma=K{\rm sph}\,e$ so that
 ${\rm tr}\,e=K^{-1}{\rm tr}\,\sigma$.
 
 Adopting the gradient theory for this internal variable $\zeta$, the 
 stored energy in terms of strain is considered
 \begin{align}\nonumber
   \mathscr{W}(u,\zeta)&=\int_\varOmega\frac12\mathbb Ce(u){:}e(u)+
   \frac12M(\beta{\rm tr}\,e(u){-}\zeta)^2
   +\frac12L(\zeta{-}\zeta_{\rm eq})^2+\frac\kappa2|\nabla\zeta|^2\,\d x
   \\&\nonumber
   =\int_\varOmega
   \frac12\Big(K+\frac{\beta^2}{d}M\Big)|{\rm sph}\,e(u)|^2+G|{\rm dev}\,e(u)|^2
   \\[-.3em]&\nonumber\qquad\qquad\qquad-\beta M\zeta{\rm tr}\,e(u)+\frac12M\zeta^2+\frac12L(\zeta{-}\zeta_{\rm eq})^2+\frac\kappa2|\nabla\zeta|^2\,\d x
 \end{align}
 which, in terms of the (here partial) stress $\sigma=\mathbb C e$,
 reads as
$\int_\varOmega\frac12(\frac1K+\frac{\beta^2}{dK^2}M)|{\rm sph}\,\sigma|^2
+\frac1G|{\rm dev}\,\sigma|^2-\beta\zeta\frac MK{\rm tr}\,\sigma+\frac12M\zeta^2
   +\frac12L(\zeta{-}\zeta_{\rm eq})^2\,\d x$.
Here $M>0$ and $\beta>0$ are so-called Biot modulus and coefficient,
 respectively, and $\zeta_{\rm eq}$ is a given equilibrium content. 
 
 We arrive at the overall stored energy as:
\begin{align}\nonumber
  \varPhi(\sigma,\zeta)=
  \int_\varOmega&\frac12\Big(\frac1K+\frac{\beta^2}{dK^2}M\Big)|{\rm sph}\,\sigma|^2
  +\frac1G|{\rm dev}\,\sigma|^2-\beta\zeta\frac MK{\rm tr}\,\sigma\,\d x
  \\[-.3em]&\qquad\qquad\qquad\qquad+
  \!\!\!\!\!\lineunder{\int_\varOmega\frac12M\zeta^2
   +\frac12L(\zeta{-}\zeta_{\rm eq})^2
  +\frac\kappa2|\nabla\zeta|^2\,\d x\,,\!}{$=:\Phi_1(\zeta)$}
\end{align}
where $\kappa>0$ is a capillarity constant.
Let us note that $\varPhi_\sigma'(\sigma,\zeta)=\bbC^{-1}\sigma+\frac{\beta M}{dK^2}(\beta{\rm sph}\,\sigma
{-}\zeta K\mathbb I)$, i.e.\ the elastic strain, and that
the proto-stress $\varSigma=\sigma$ indeed differs from an actual stress
by the spherical pressure part $\frac{\beta M}{dK}(\beta{\rm sph}\,\sigma
{-}\zeta K\mathbb I)$.\COMMENT{HERE THE CALCULATIONS BETTER TO CHECK}

The driving force for the diffusion is the
chemical potential $\mu=\varPhi_\zeta'(\sigma,\zeta)$, i.e.\ here
\begin{subequations}\label{CH}\begin{align}
  \mu=(M+L)\zeta-\beta \frac MK{\rm tr}\,\sigma
  -L\zeta_{\rm eq}-\kappa\Delta\zeta\,.
\end{align}
The diffusion equation is
\begin{align}
  \DT\zeta-{\rm div}(\mathbb M\nabla\mu)=0
\end{align}\end{subequations}
with $\mathbb M$ denoting the diffusivity tensor.
The system \eqref{CH} is called the Cahn-Hilliard equation, here combined
with elasticity so that the flow of the diffusant is driven both
by the gradient of concentration (Fick's law) and the gradient of
the mechanical pressure (Darcy's law). The dissipation potential in terms of
$\nabla\mu$, let us denote it by $R$ behind this system is
\begin{align}
R(\mu)=\int_\varOmega\frac12\mathbb M\nabla\mu{\cdot}\nabla\mu\,\d x,
\end{align}
For the analysis cf.\ e.g.\ \cite[Sect.\,7.6]{KruRou18MMCM}.

One would expect the dissipation potential as a function of the rate
of internal variables, as in \eqref{IVP2z+}. In fact,
the system \eqref{CH} turns into the form \eqref{IVP2z+} is one takes
the dissipation potential $\varPsi=\varPsi(\DT\zeta)$ as
\begin{align}
\varPsi(\DT\zeta)=R^*(\DT\zeta)
\end{align}
with $R^*$ denoting the convex conjugate to $R$. Now,
$\varPsi$ is nonlocal. 
The functional setting is as in the previous example but now
$\calZ=H^1(\varOmega)$ and $\calZ_1=H^1(\varOmega)^*$.
For a discretisation of the
type \eqref{suggestion-3}, see \cite{Roub17ECTD}.

Often, the diffusivity is considered dependent on $\zeta$.
Or even one can think about $\mathbb M=\mathbb M(\sigma,\zeta)$. Then the
modification in Remark~\ref{rem-dissip} is to be applied. In particular,
$R(\sigma,\zeta,\mu)=\int_\varOmega\frac12\mathbb M(\sigma,\zeta)\nabla\mu
  {\cdot}\nabla\mu\,\d x$ and $\varPsi(\sigma,\zeta,\DT\zeta)=[R(\sigma,\zeta,\cdot)]^*(\DT\zeta)$.

  For this Biot model in the dynamical variant, the reader is also referred to
  the books \cite{AbChUl05PBC,Carc15WFRM,Chen16P,Stra08SWMP} or also
  \cite{KruRou18MMCM,MieRou15RIST}. In any case, the diffusion involves
  gradients and in the implicit discretisation it leads to large systems of
  algebraic equations, which inevitably slows down the fast explicit
  discretisation of the mechanical part itself.

 \subsection{Damage}

The simplest examples of nonconvex stored energy are models of damage.
The most typical models use as an 
internal variable the scalar-valued bulk damage $\alpha$ having the
interpretation as a phenomenological 
volume fraction of microcracks or microvoids manifested macroscopically
as a certain weakening of the elastic response.
This concept was invented by L.M.\,Kachanov \cite{Kach58TRPD}
and Yu.N.\,Rabotnov \cite{Rabo69CPSM}.

Considering gradient theories, the stored energy in terms of the strain and
damage is here considered as
$$
\mathscr W(e,\alpha)
=\int_\varOmega\frac12\gamma(\alpha)\mathbb C e{:}e+\phi(\alpha)
+\frac\kappa2|\nabla\alpha|^2+\frac\varepsilon2\nabla(\mathbb C e)
{:}\nabla e\,\d x\,,
$$
where $\phi(\cdot)$ is an energy of damage which gives rise to an activation
threshold for damage evolution and may also lead to healing (if allowed).
The last term is mainly to facilitate the mathematics towards convergence
and existence of a weak solution in such purely elastic materials
without involving any viscosity, cf.\
\cite[Sect.\,7.5.3]{KruRou18MMCM}. This regularization can also control
dispersion of elastic waves.
The $\nabla\alpha$-term also facilitates the analysis and controls
the internal length-scale of damage profiles.

Let us consider the ``generalized'' elasticity tensor $\mathfrak C=\bbC$
independent of $x$. As in the previous examples, $Eu=e(u)$ and $G=0$.
According \eqref{IVP2Sigma+}, the proto-stress $\varSigma=\mathfrak CEu+G$,
denoted by $\sigma$, now looks as $\mathbb C e=:\sigma$; in damage
  mechanics, the proto-stress is also called an effective stress with
  a specific mechanical interpretation, cf.\ \cite{Rabo69CPSM}.
In terms of $\sigma$, the stored energy is then
  \begin{align}
  \varPhi(\sigma,\alpha)=
  &\int_\varOmega\frac12\gamma(\alpha)\mathbb C^{-1}\sigma{:}\sigma
  +\frac\varepsilon2\nabla\mathbb C^{-1}\sigma{:}\nabla\sigma\,\d x
 +\!\!\!\!\lineunder{\int_\varOmega
  \phi(\alpha)+\frac\kappa2|\nabla\alpha|^2\,\d x\,.\!}{$=:\Phi_1(\alpha)$}
  \label{damage-stored}
  \end{align}
Then $\varPhi_\sigma'=\gamma(\alpha)\bbC^{-1}\sigma
    -{\rm div}(\varepsilon\nabla(\bbC^{-1}\sigma))$ and
    the true stress $S=\bbC^*\varPhi_\sigma'$ is then
    $\gamma(\alpha)\sigma-{\rm div}(\varepsilon\nabla\sigma)$ provided $\bbC$ is
    constant and symmetric.
 The damage driving force (energy) is $\varPhi_\alpha'(\sigma,\alpha)=
    \frac12\gamma'(\alpha)\mathbb C^{-1}\sigma{:}\sigma+\phi'(\alpha)
    -{\rm div}(\kappa\nabla\alpha)$. When $\gamma'(0)=0$ and $\phi'(0)\le0$,
    then always $\alpha\ge0$ also in the discrete scheme if $\alpha_0\ge0$.

    The other ingredient is the dissipation potential. To comply with the
    coercivity on $\calZ_1=L^2(\Omega)$ with $p\ge2$ as needed in
    Proposition~\ref{prop1}, one can consider either
 \begin{align}\label{Psi-damage}  
   \varPsi(\DT\alpha)=\begin{cases}\int_\varOmega\varepsilon_1\DT\alpha^2\,\d x&
     \\\ \ +\infty&\end{cases}
     \text{ or }\ \ \ \begin{cases}\int_\varOmega\varepsilon_1\DT\alpha^2\,\d x&\text{if }\DT\alpha\le0\text{ a.e.\ on }\varOmega,
     \\\int_\varOmega\DT\alpha^2/\varepsilon_1\,\d x&\text{otherwise}\end{cases}
 \end{align}
 with some (presumably small) coefficient $\varepsilon_1>0$. The former option
 corresponds to a unidirectional (i.e.\ irreversible) damage not allowing any
 healing (as used in engineering) while the latter option allows for
 (presumably slow) healing as used in geophysical models on large time scales.
 
    Since $\sigma$ appears nonlinearly in $\varPhi_\alpha'(\sigma,\alpha)$,
the strong convergence $\overline\sigma_{\tau h}\to\sigma$ in $L^2(Q;\R^{d\times d})$
    is needed. For this, the strain-gradient term with $\varepsilon>0$ is needed
    and the Aubin-Lions compact embedding theorem is used. This gives
    the strong convergence even in the norm of
    $L^{1/\epsilon}(I;L^{2d/(d-2)-\epsilon}(\varOmega;\R^{d\times d}))$ for
    arbitrarily small $\epsilon>0$ provided also $\DT\sigma_{\tau h}$ is bounded
    in some norm, which can be shown by using
    $\DT\sigma_{\tau h}=\bbC e(\overline v_{\tau h})$ and the Green formula
    \begin{align}\nonumber
      \big\|\DT\sigma_{\tau h}\|_{L^\infty(I;H^{-1}(\varOmega;\R^{d\times d}))}^{}
      &=\sup_{\|\widetilde e\|_{L^1(I;H_0^1(\varOmega;\R^{d\times d}))}\le1}
      \int_0^T\!\!\int_\varOmega\DT\sigma_{\tau h}{:}\widetilde e\,\d x\d t
      \\&\nonumber=\sup_{\|\widetilde e\|_{L^1(I;H_0^1(\varOmega;\R^{d\times d}))}\le1}
     \int_0^T\!\!\int_\varOmega\bbC e(\overline v_{\tau h}){:}\widetilde e\,\d x\d t
      \\&\nonumber=\sup_{\|\widetilde e\|_{L^1(I;H_0^1(\varOmega;\R^{d\times d}))}\le1}\!\!
-\int_0^T\!\!\int_\varOmega\overline v_{\tau h}{\cdot}{\rm div}(\bbC\widetilde e)\,\d x\d t
      \le C\|\overline v_{\tau h}\|_{L^\infty(I;L^2(\varOmega;\R^d))}^{}
    \end{align}
    with $C$ depending on $|\bbC|$. Cf.\ also the abstract estimation
    \eqref{dual-est1}.

When $\gamma$ or $\phi$ are not quadratic but continuously
differentiable, one can use the abstract difference quotient \eqref{quotient}
defined, in the classical form, as
\begin{align}
  \varPhi_z^\circ(\varSigma,\alpha,\widetilde\alpha)=
  \begin{cases}\frac12\frac{\gamma(\alpha){-}\gamma(\widetilde\alpha)}{\alpha{-}(\widetilde\alpha)}
    \mathbb C^{-1}\sigma{:}\sigma
    +\frac{\phi(\alpha){-}\phi(\widetilde\alpha)}{\alpha{-}(\widetilde\alpha)}
    -\kappa\Delta\frac{\alpha{+}\widetilde\alpha}2
      &\text{where }\ \alpha\ne\widetilde\alpha\,.
      \\\frac12\gamma'(\alpha)\mathbb C^{-1}\sigma{:}\sigma+\phi'(\alpha)
        -\kappa\Delta\alpha&\text{where }\ \alpha=\widetilde\alpha\,.
    \end{cases}\label{nonquadratic-damaga}
\end{align}
Of course, rigorously, the $\Delta$-operator in \eqref{nonquadratic-damaga}
is to be understood in the weak form when using it in \eqref{suggestion-3}.

Due to the gradient $\kappa$-term in \eqref{damage-stored}, the implicit
incremental problem \eqref{suggestion-3} leads to an algebraic problem
with a full matrix, which may substantially slow down the otherwise fast
explicit scheme. Like in the previous model the capillarity, now this gradient
theory controls the length-scale of the damage profile and also serves as
a regularization to facilitate mathematical analysis. Sometimes, a nonlocal
``fractional'' gradient
can facilitate the analysis, too.
Then, some
wavelet equivalent norm can be considered to accelerate the calculations,
cf.\ also \cite{ArGrRo03MNSM}.
  As far as the stress-gradient term,
  it is important that the discretisation of the proto-stress
  in the usual implementation of the leap-frog method is continuous piecewise
  smooth, so that $\nabla\sigma$ has a good sense in the discretisation
  without need to use higher-order elements.
    Here we use that the latter relation in \eqref{suggestion-2} is to be
    understood in the weak form, namely
  $\int_\varOmega S_{\tau h}^{\KK}{:}\tilde E_h\,\d x
  =\langle\varPhi_\varSigma'(\varSigma_{\tau h}^{\KK},z_{\tau h}^{\KKK}),
  \mathfrak C\tilde E_h\rangle$ for
  $\tilde E_h=\tilde e_h=e(\tilde u_h)$, which means
  $$
  \int_\varOmega S_{\tau h}^{\KK}{:}\tilde E_h\,\d x
  =\int_\varOmega\gamma(\alpha_{\tau h}^k)\mathbb C^{-1}\sigma^{\KK}{:}
  \mathbb C\tilde e_h
  +\epsilon\nabla\mathbb C^{-1}\sigma^{\KK}{:}\nabla\mathbb C\tilde e_h\,\d x
  $$
  for any $\tilde e_h$ from the corresponding finite-dimensional subspace of
  $H^1(\varOmega;\R_{\rm sym}^{d\times d})$.
  Thus we indeed do not need higher-order elements, and also we do not need
  to specify explicitly homogeneous boundary conditions in this boundary-value problem.

  The functional setting is $\calH=L^2(\varOmega;\R^d)$, $\calE=\calS
  =H^1(\varOmega;\R_{\rm sym}^{d\times d})$,
  $\calZ=H^1(\varOmega)$, and $\calZ_0=\calZ_1=L^2(\varOmega)$.
  Then $\calU=H^2(\varOmega;\R^d)$, and $E=e(\cdot)$
  is understood as an operator $H^2(\varOmega;\R^d)\to
  H^1(\varOmega;\R_{\rm sym}^{d\times d})$, and $\mathfrak C^*\cong\bbC^\top\!=\bbC$
  is understood as an a operator from $H^1(\varOmega;\R_{\rm sym}^{d\times d})$
  to itself.


  \def\Frakg{{\mathfrak{g}}}
  A particular case of this model is a so-called {\it phase-field fracture},
  taking as basic choice
  \begin{align}\gamma(\alpha):={\varepsilon^2}/{\varepsilon_0^2}
        {+}\alpha^2,\ \ \ \
        \phi(\alpha):=\Frakg_{\rm c}{(1{-}\alpha)^2}/{\varepsilon},\ \
        \text{ and }\ \ \kappa:=\varepsilon \Frakg_{\rm c}
        \label{AT}\end{align}
  with $\Frakg_{\rm c}$ denoting the energy of fracture and with $\varepsilon$
  controlling a ``characteristic'' width of the {\it phase-field fracture}.
  The physical dimension of $\varepsilon_0$ as well as of 
$\varepsilon$ is m (meters) while the physical dimension of $\Frakg_{\rm c}$ 
is J/m$^2$. This is known as the so-called \emph{Ambrosio-Tortorelli 
  functional} \cite{AmbTor90AFDJ}. In the dynamical context, only various
implicit discretisation schemes seems to be used so far, cf.\
\cite{BVSH12,HofMie12CPFM,RouVod??MMPF,SWKM14PFAD}.
There are a lot of improvements of this basic model, approximating 
a mode-sensitive fracture, or $\varepsilon$-insensitive models
(with $\varepsilon$ referring to \eqref{AT}), or ductile fracture. This last
variant combines this model with the plasticity as in Sect.~\ref{sect-plast}.

  \subsection{Delamination on adhesive contacts}\label{sec-delam}

  Let us now present an example for a less trivial operator $E$, namely
  $Eu=(e,w)$ with $e=e(u)$ on $\varOmega$ as before and with $w=u|_\varGamma)$
  being the trace of $u$ on the boundary $\varGamma$. 
  The internal variable will be a scalar-valued surface damage $\alpha$ on
  $\varGamma$, i.e.\ the so-called delamination variable, which is the
  concept introduced by M.\,Fr\'emond \cite{Frem85DLS}.
  
 The stored energy in terms of the strain and trace of the displacement is
 $$
\mathscr W(e,w,\alpha)=\int_\varOmega\frac12\mathbb C e{:}e\,\d x
 +\int_\varGamma\frac12\gamma(\alpha)\bbB
w{\cdot}w+\COLOR{\phi(\alpha)+}\frac\kappa2|\nablaS\alpha|^2\,\d S$$
 with $\bbC\in\R^{d\times d\times d\times d}$ 
 symmetric positive definite and $\bbB\in\R^{d\times d}$
 symmetric positive semidefinite, and with
 $\nablaS$ a surface gradient. This leads us to
 consider $\mathfrak C=\mathbb C\times\bbB$ and the proto-stress
 $\varSigma:=(\sigma,\varsigma)$.
 The stored energy expressed in terms of this proto-stress is
\begin{align}
  \varPhi(\sigma,\varsigma,\alpha)=
  \int_\varOmega\frac12\mathbb C^{-1}\sigma{:}\sigma\,\d x
+\int_\varGamma\frac12\gamma(\alpha)\mathbb\bbB^{-1}\varsigma{\cdot}\varsigma\,\d S
  +\!\!\!\!\lineunder{\int_\varGamma\COLOR{\phi(\alpha)+}\frac\kappa2|\nablaS\alpha|^2\,\d S\,.\!}
  {$=:\Phi_1(\alpha)$}
\label{Phi-delam}\end{align}
The dissipation potential is usually taken as in \eqref{Psi-damage} except
that $\varOmega$ is replaced by $\varGamma$. The damage gradient term in
\eqref{Phi-delam} leads to the  Laplace-Beltrami operator on the boundary in
the classical formulation of the flow-rule for $\alpha$.

Moreover we may consider the boundary
loading through the Robin boundary condition
$\sigma\vec{n}=\gamma(\alpha)\bbB(u_{\rm D}^{}(t){-}u)$ on $\varGamma$ with some given
displacement $u_{\rm D}^{}$ depending on time. Then $G(t)\in\calS$ is given by
$\langle G(t),(\sigma,\varsigma)\rangle
=-\int_\varGamma \bbB u_{\rm D}^{}\cdot\varsigma\,\d S$.
The bulk load $F\in\calH^*$ is considered as
$\langle F,u\rangle=\int_\varOmega f\cdot u\,\d x$.

Thus $\varPhi_\varSigma'(\sigma,\varsigma,\alpha)
=(\mathbb C^{-1}\sigma,\gamma(\alpha)\bbB^{-1}\varsigma)
$ and
the generalized actual stress is $S=\mathfrak C^*\varPhi_\varSigma'=
(\sigma,\gamma(\alpha)\varsigma)
$.
The abstract identity $\varSigma={\mathfrak C}Eu+G$ occurring in
Definition~\ref{def} means component wise that
\begin{align}
  \sigma=\bbC e(u)\ \ \ \text{ and }\ \ \
  \varsigma=\bbB(u_{\rm D}^{}-u|_\varGamma).
 \label{delam-sigma}\end{align}
The abstract force equilibrium $\mathcal T'\DT v+E^*S 
=F(t)$, cf.\ \eqref{IVP2S+},  with the initial condition $v(0)=v_0$
in the weak form \eqref{weak-v} gives 
\begin{align}
  \int_0^T\!\!\int_\varOmega\sigma:e(\widetilde v)-\varrho v\cdot\DT{\widetilde v}
  \,\d x\d t+\int_0^T\!\!\int_\varGamma\gamma(\alpha)\varsigma\cdot\widetilde v\,\d S\d t
  =\int_\varOmega\varrho v_0\cdot\widetilde v(0)\,\d x
  +\int_0^T\!\!\int_\varOmega f\cdot \widetilde v\,\d x\d t.
  \label{delam-equil}\end{align}
Substituting \eqref{delam-sigma} into \eqref{delam-equil} and
taking into account that $v=\DT u$, we obtain
the weak formulation of the equation \eqref{IBVP1} with the
initial condition \eqref{IBVP2} and the 
boundary condition $[\bbC e(u)]\vec{n}+\gamma(\alpha)\bbB u=\gamma(\alpha) u_{\rm D}^{}$.
In particular, we can also see that $\sigma\vec{n}=\varsigma$.

Here the functional setting is
$\calH=L^2(\varOmega;\R^d)$, $\calU=H^1(\varOmega;\R^d)$,
$\calE=\calS=L^2(\varOmega;\R_{\rm sym}^{d\times d})\times L^2(\varGamma;\R^d)$,
$\calZ=H^1(\varGamma)$, and $\calZ_0=\calZ_1=L^2(\varGamma)$.


When the adhesive is close to be brittle (i.e.\ $\bbB$ is big), the
CFL-condition \eqref{CFLnew} becomes very restrictive.  
For a ``stabilization'' of the explicit method for the such brittle adhesive,
one can use an artificial mass on the boundary, cf.\
\cite{KTGKBP17SRTA,PaPaMa10ITDB,PaPaMa11CAPT}. This spurious mass can,
however, be suppressed to zero if the CFL condition is strengthened so that
$\tau/h\to0$.

Let us still note the that Neumann boundary conditions
  can easily be considered instead of the Robin boundary conditions. 
Also the surface-gradient term in \eqref{Phi-delam} can
be omitted if both $\gamma(\cdot)$ and $\phi(\cdot)$ are affine, the latter one
being still augmented by the indicator function of the interval $[0,1]$
to ensure that $\alpha$ is valued in this interval, cf.\ e.g.\
\cite{MieRou15RIST,RKVPZ??DACM}. Then $\calZ=L^2(\varGamma)$ and
the latter option in (\ref{ass-Phi'}a,b) is to be used, and
even the equation (or inclusion)
\eqref{suggestion-3} is local (like in Sect.\,\ref{sect-plast} without
gradient of plastic strain) and the discretisation is the truly explicit.
This will be used in Sect.\,\ref{sec-numer}.

The model presented so far has limited application because, after a complete
delamination of the adhesive contact, such part of the boundary becomes
completely free and allows unrealistically for the penetration with the
obstacle. A simple improvement of this model combines the damageable Robin
boundary condition in the tangential direction with homogeneous Dirichlet
boundary condition in the normal direction. This leads to a so-called
{\it bilateral contact}.

\section{Implementation and 2D-numerical experiments}\label{sec-numer}

In this last section, we demonstrate the efficiency of the explicit discretisation (which is well recognized for the linear elastodynamics) and combined with dissipative evolution of internal variables in the staggered way, as devised above. We use the delamination model from Sect.\,\ref{sec-delam}.

For the discretisation, we use the lowest order $Q^{\rm div}_{k+1}-Q_k$ finite element (for $k=0$) proposed and analyzed in \cite{BeJoTs02NFMF} for the linear elastodynamic problem written as a first order hyperbolic system with unknowns the velocity $v$ and the stress tensor $\sigma$. We consider the two dimensional problem and the domain is discretised with rectangular elements.  For the stress tensor we use piecewise bi-linear functions with the following continuities: the normal stress is continuous across edges of adjacent elements while the tangential component of the stress may be discontinuous. The velocity is discretised by piecewise constant functions. Similarly, considering $\kappa=0$ in \eqref{Phi-delam}, 
  we use the $P_0$-elements on the side segments for discretisation of the delamination variable $\alpha$.
For more details about the elastodynamic part,  we refer to \cite{BeJoTs02NFMF} and \cite{Tsog99}. The extension here is that Robin-type boundary conditions have been incorporated to the scheme that are quite general and may serve to appropriately describe mixed (stress-velocity) conditions on the boundary as well as the prescribed stress and velocity conditions.

In the following example we present the delamination model from Sect.\,\ref{sec-delam} on an adhesive contact case.
\begin{figure}
       \centerline{\hspace*{0em}
 \includegraphics[width=.8\textwidth]{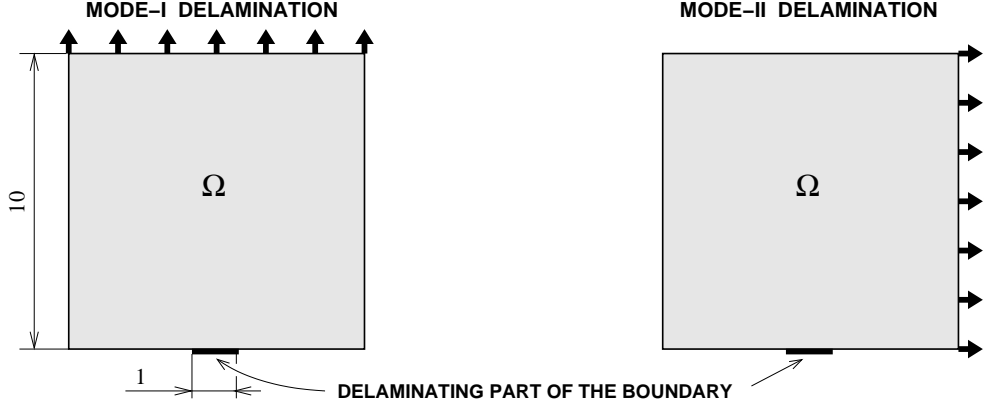}}
       \caption{{\sl
          The square-shaped 2-dimensional domain
$\varOmega$ and the boundary conditions:
the mid-part of the bottom side
           where damageable Robin boundary conditions hold
           (i.e.\ the adhesive contact),
           traction boundary conditions are considered in the rest,
           either homogeneous or gradually increasing in time, considering
           two options
           leading primarily to Mode I and Mode II as depicted in the left
       and the right figure, respectively.}}
       \label{fig:compdom+}
\end{figure}
For mere demonstration of efficiency of the time/space discretisation
and the algorithm, we take dimensionless data, i.e.\ without physical units. 
The material of the specimen considered here is assumed to be homogeneous and isotropic, with
the bulk modulus $K{=}$1.66 and the shear modulus $G{=}$1.
  We further consider $\varrho=1$, which then corresponds to the pressure-wave velocity $v_{\rm p}=\sqrt{(K{+}4G/3)/\varrho}=\sqrt3=1.73$ and
  the shear-wave velocity $v_{\rm s}=\sqrt{G/\varrho}=1$.

The size of the square-shaped domain is as depicted in Figure~\ref{fig:compdom+}. Namely, a square $\varOmega=(0,10)^2$ and
the part of the boundary amenable to damage and thus debonding from
an outer support 
is located at the center of the bottom
boundary, cf.\ Fig.~\ref{fig:compdom+}. Boundary
conditions are of the form \eqref{IBVP2+} with $\DT g{=}0$ and the adhesive stiffness $\bbB =\frac12\bbI$.

Furthermore, we took $\gamma(\alpha)=\alpha$ and considered the explicit
  constraints $0\le\alpha\le1$, as mentioned in Sect.~\ref{sec-delam} as an
  alternative.
The dissipation potential $\Psi$ of~\eqref{IVP2z+} is
taken like the former option in \eqref{Psi-damage} with
  $\varepsilon_1=0$ and $\varOmega$ replaced by the delaminating
  part of $\varGamma$, and $\phi(\alpha)=\mathfrak{g}\alpha$
  with the fracture-toughness constant $\mathfrak{g}{=}2.57\cdot10^{-5}$.

  The excitation is imposed on the loaded side by normal stresses assumed to vary linearly
  and slowly in time as $\sigma_n{=}0.005t/T$, with $T{=}51$ being the total duration of
  the experiment. The tangent traction is assumed to be zero while, at the rest of the boundary,
  \COLOR{the 
    traction-free} boundary conditions are assumed and enforced. The computational experiments have
  been performed with the mesh size $h{=}0.025$, that gives a grid of $400\times 400$ elements.
  Thus, on adhesive-contact boundary, there are $40$ elements. The time discretisation step is $\Delta t= 
\frac{1}{v_p} h =0.0144$. 

The standard Helmholtz decomposition \cite{achenbach2012wave}, usually used (e.g. in seismology~\cite{aki2002quantitative}), is performed. The computed velocity gradient is therefore decomposed in its pressure and shear waves using $\nabla v=({\rm div}\,v){\mathbb I}/d+{\rm rot}\,v$.

\begin{figure}
  \ \\[-2.2em]
   \includegraphics[trim=100 200 160 100,clip, width=.8\textwidth]{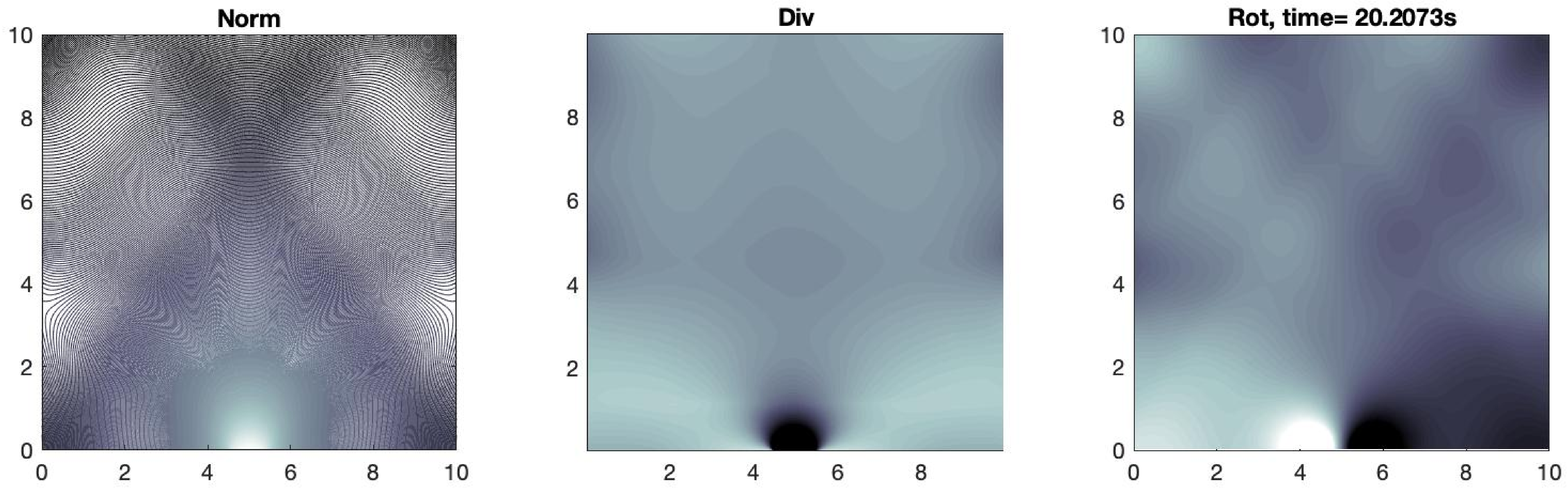}\\[-1.2em]
 \includegraphics[trim=100 200 160 100,clip,width=.8\textwidth]{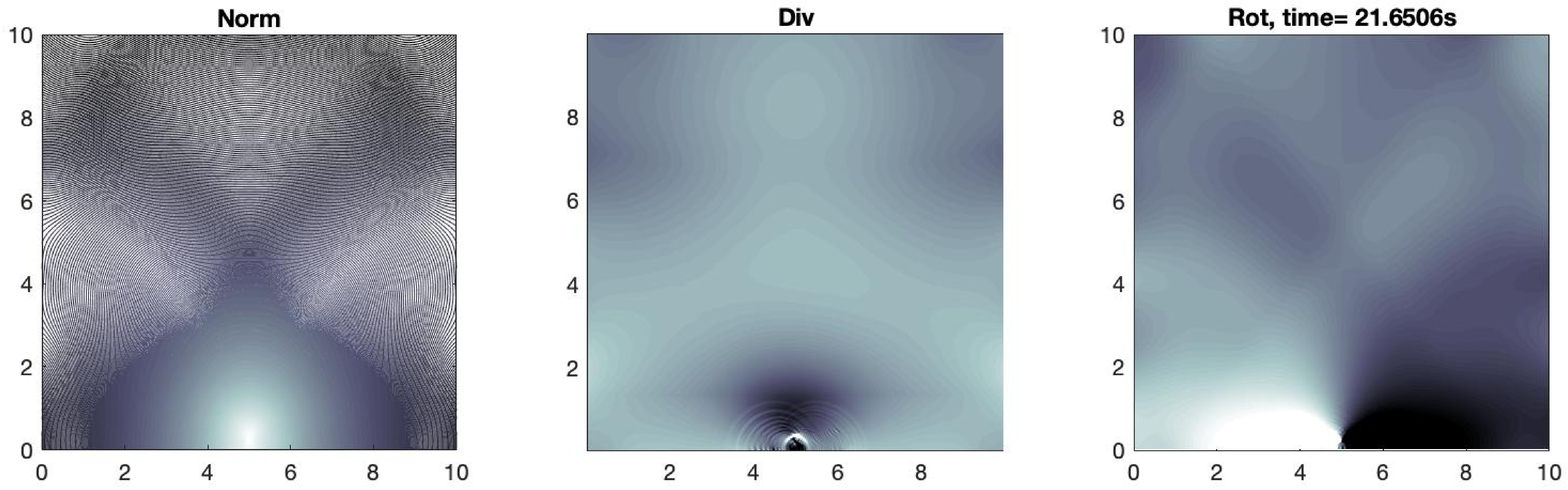}\\[-1.2em]
 \includegraphics[trim=100 200 160 100,clip,width=.8\textwidth]{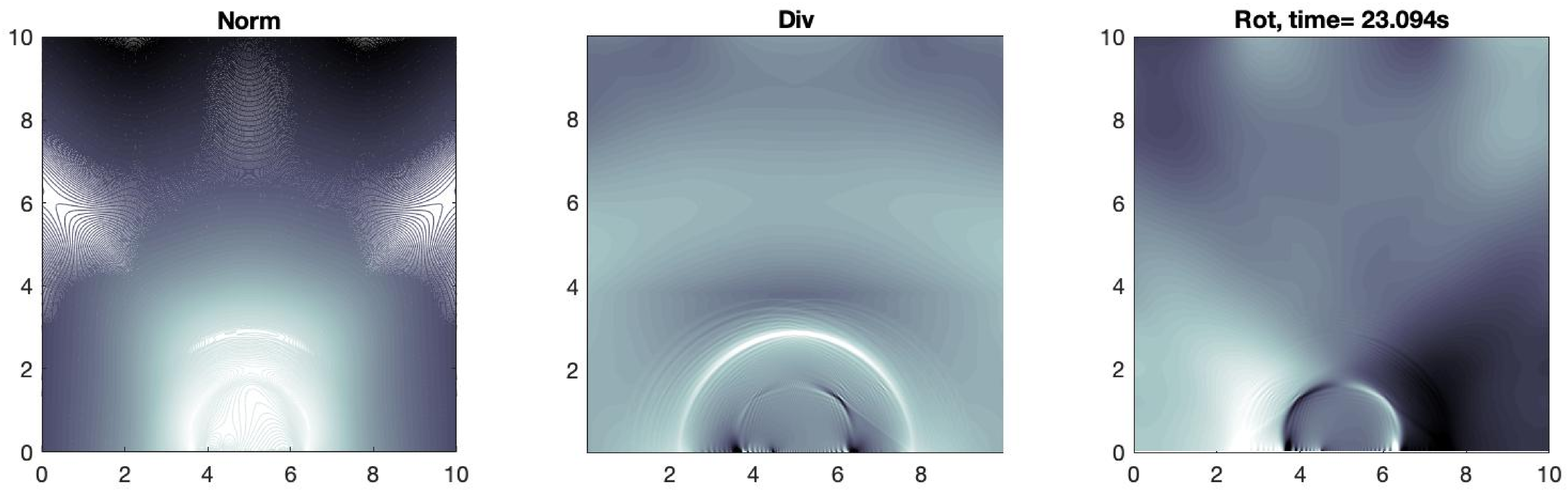}\\[-1.2em]
 \includegraphics[trim=100 200 160 100,clip,width=.8\textwidth]{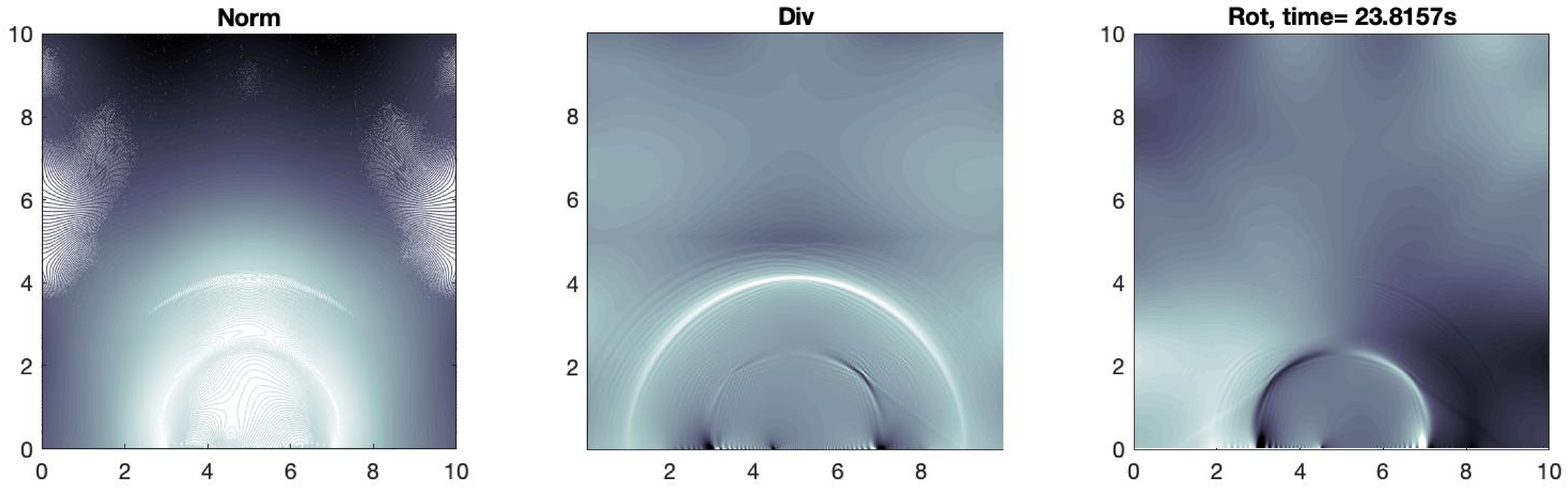}\\[-1.2em]
 \includegraphics[trim=100 200 160 100,clip,width=.8\textwidth]{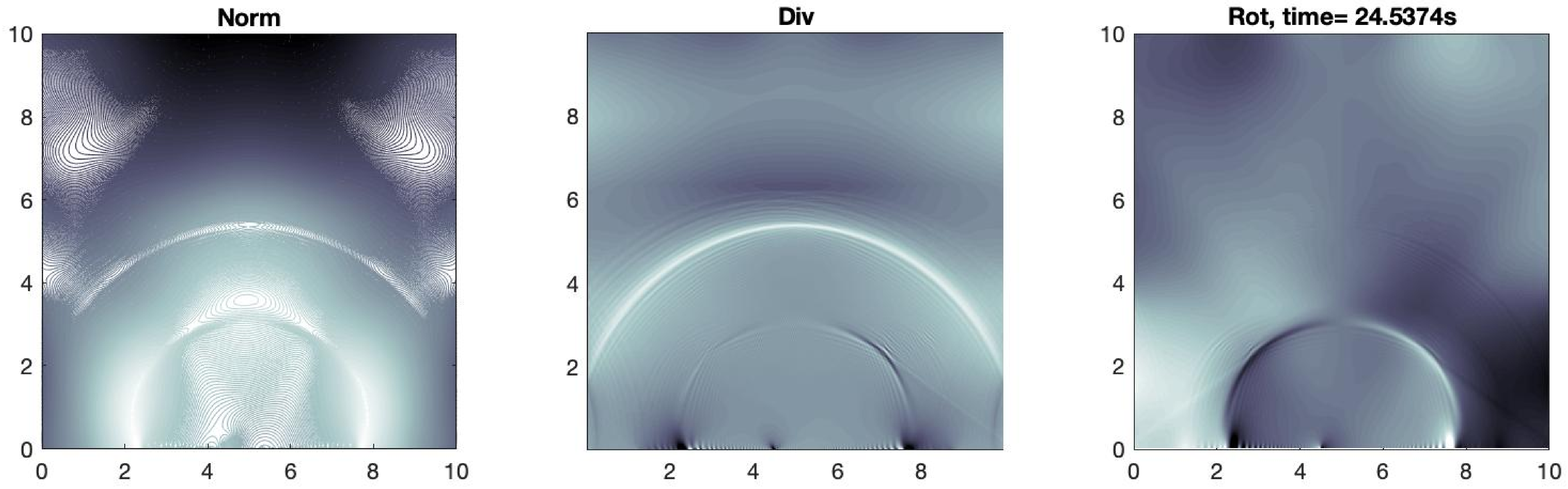}\\[-2.5em]
 \caption{{\sl Delamination (rather) in Mode I under the loading as in
     Fig.\,\ref{fig:compdom+}-left.
         Five selected time instants
   immediately after the delamination was executed are displayed in the following rows.
   Each row consist in spatial distribution of the norm of velocity,
   divergence of velocity, and rotation of velocity. 
   Both P-wave and S-waves are
   emitted, the former one being faster, as clearly seen on the middle column.}
   \label{fig:test-1}}
\end{figure}

In the first experiment, the top side of the square-shaped domain is loaded as shown in Fig.~\ref{fig:compdom+}(left). As it was expected a rather Mode-I debonding takes place.  The acoustic emission generated by fast propagation of surface damage can be seen in Fig.~\ref{fig:test-1} where the norm of the velocity vector field is plotted together with both the divergence and its rotational part to identify P-wave and S-wave, respectively. The term ``acoustic emission'' is used to describe the transient elastic waves caused by the rapid release of localized stress energy, but one can also understand it as a seismic-wave emission, depending on a particular application. This localization of the stress energy can be seen on the top row plots of Fig.~\ref{fig:test-1} corresponding  to time $t=20.2073$. Rupture is occurring rather rapidly during a very short time of successive symmetrical appeared damage events. Then elastic waves (both pressure and shear) emanate from the damaged region and propagate \REPLACE{inside the medium}{thorough the specimen}
as illustrated on selected snapshots at times $t=21.6506$, $t=23.094$, $t=23.8197$, and \COLOR{$t=24.630$} in Fig.~\ref{fig:test-1}. 

\begin{figure}
 \includegraphics[trim=100 200 160 100,clip,width=.8\textwidth]{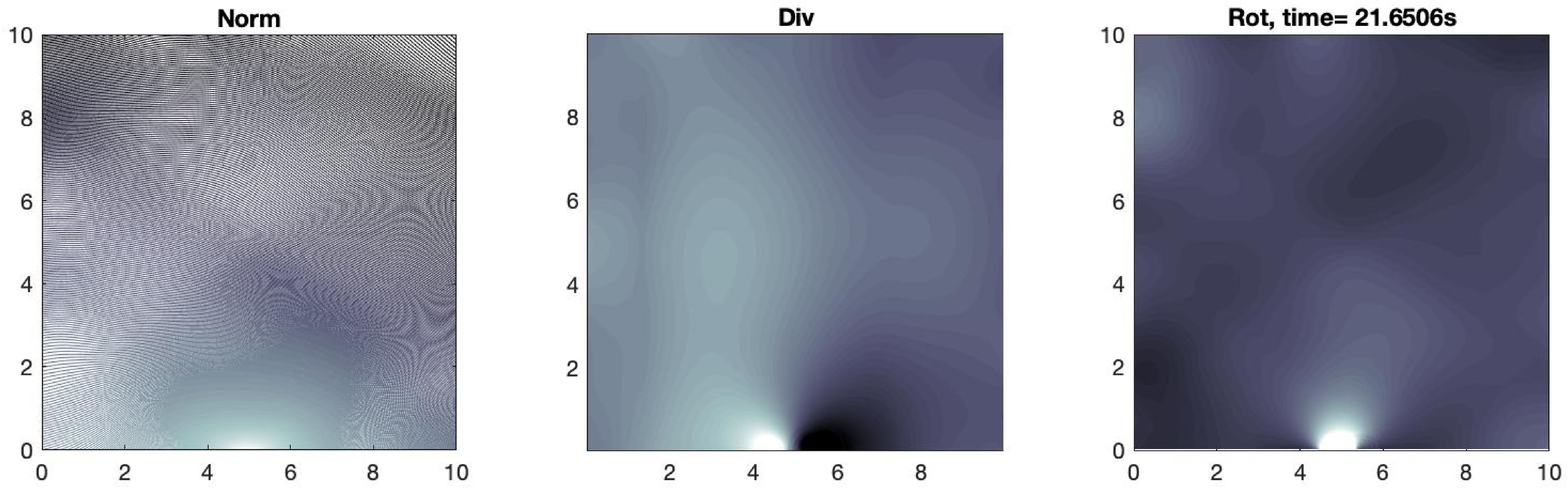}\\[-1.2em]
 \includegraphics[trim=100 200 160 100,clip,width=.8\textwidth]{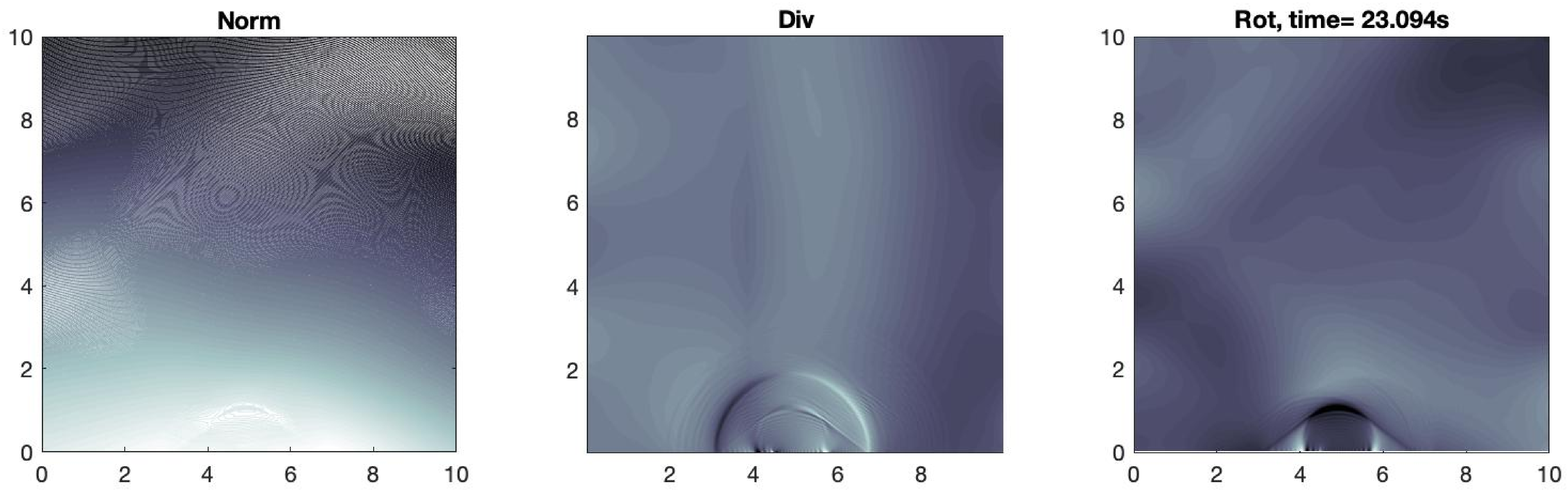}\\[-1.2em]
 \includegraphics[trim=100 200 160 100,clip,width=.8\textwidth]{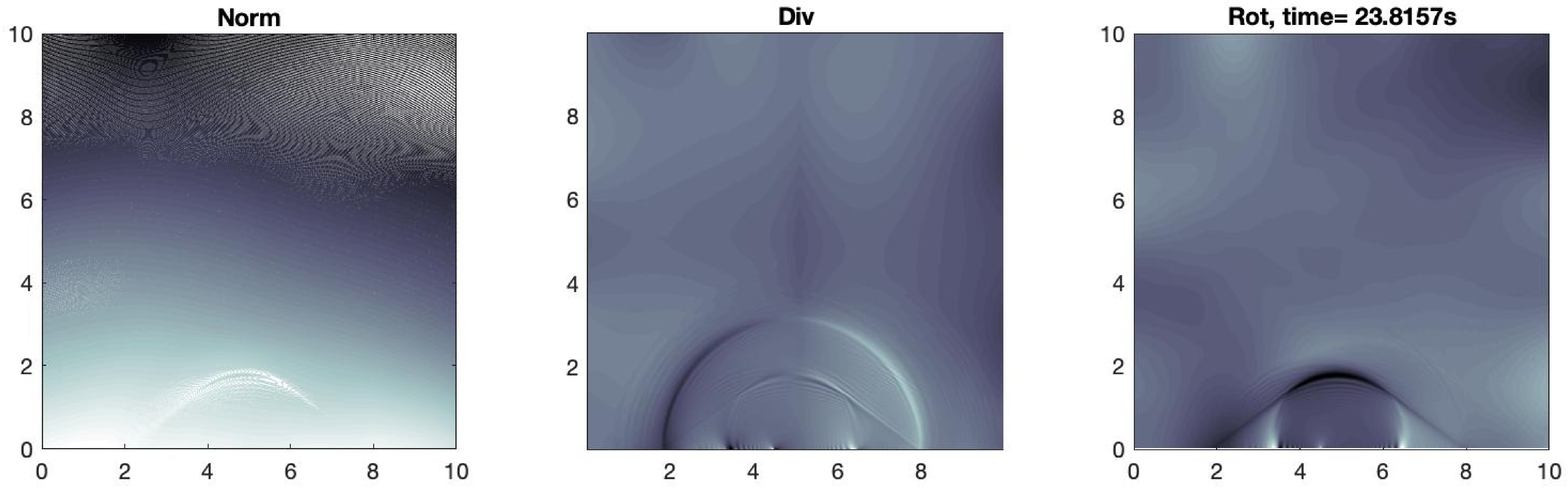}\\[-1.2em]
 \includegraphics[trim=100 200 160 100,clip,width=.8\textwidth]{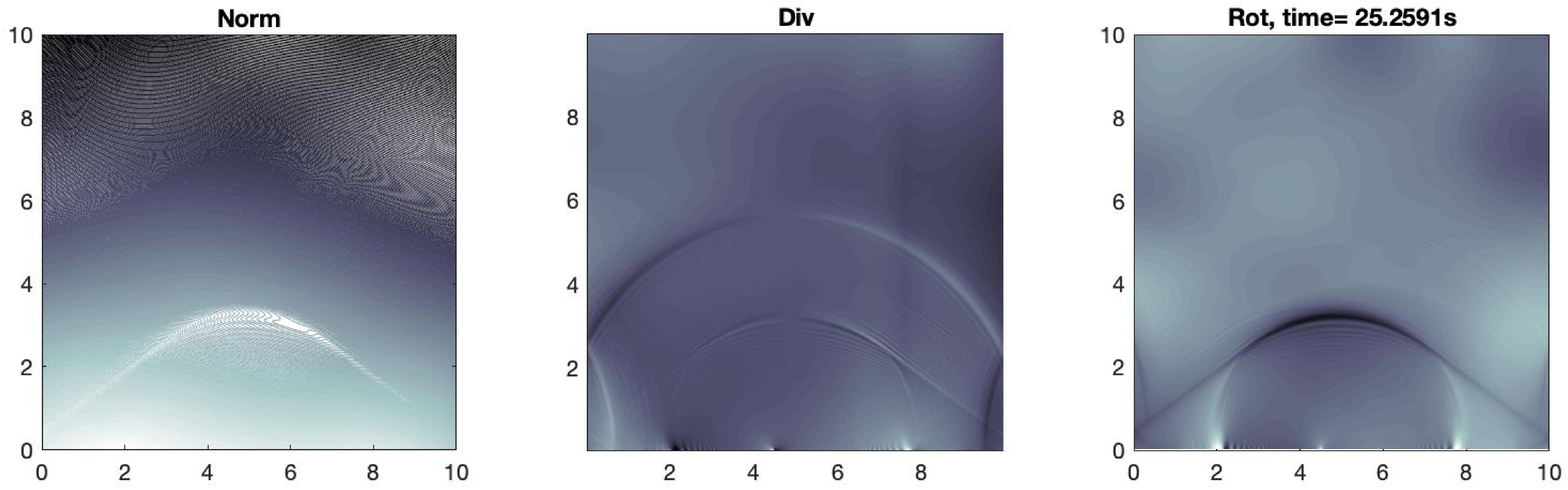}\\[-1.2em]
 \includegraphics[trim=100 200 160 100,clip,width=.8\textwidth]{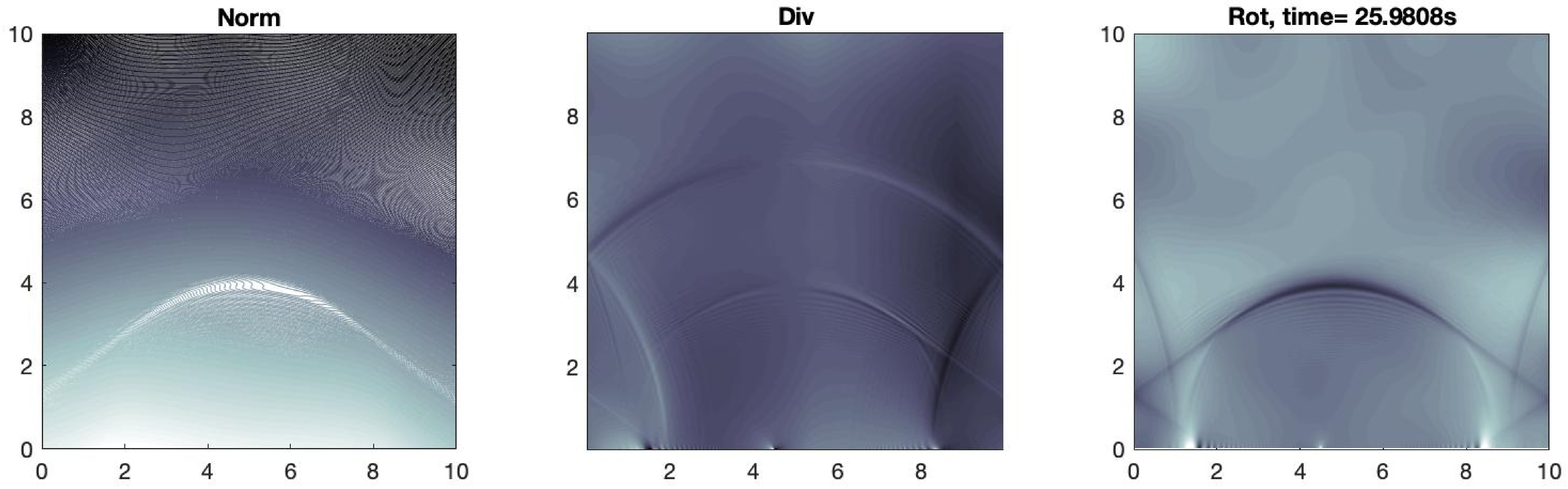}\\[-2.5em]
 \caption{{\sl Delamination (rather) in Mode II under the loading as in
   Fig.\,\ref{fig:compdom+}-right. In the left column, it is
   clearly visible that the S-wave dominates.} \label{fig:test-2}}
\end{figure}

In the second experiment, a rather Mode-II damage evolution is performed by imposing the loading pattern of Fig.~\ref{fig:compdom+}(right). In that case, a bit longer duration of time is needed for the damage to start expanding on  the adhesive part of the boundary. The localization of the stress energy is illustrated on the top row plots of Fig.~\ref{fig:test-2} corresponding  to time $t=21.6506$. The rupture in this case consists of non symmetrical occurrences of damage events. The wave propagation of the velocity field can be seen in the plots of Fig.~\ref{fig:test-2} with a more articulated S-wave while $P$-wave is rather suppressed.

In both experiments, the waves are quite sharp and nearly without spurious dispersion during their propagation. This shows efficiency of the explicit numerical algorithm even if combined with the dissipative inelastic process.  

%


\bigskip\bigskip

\baselineskip=12pt

\noindent{\small
{\it Acknowledgments.}
\COLOR{Fruitful discussions with Dr.~Radek Kolman are warmly acknowledged.}
This research has been partly supported from the grants 17-04301S (especially
as far as the focus on the dissipative evolution of internal variables)
and 19-04956S (especially as far as the focus on the dynamic and nonlinear
behaviour)
of the Czech Sci.\ Foundation. C.G.P.\ acknowledges the financial support of
the Stavros Niarchos Foundation within the framework of the project ARCHERS
(``Advancing Young Researchers' Human Capital in Cutting Edge Technologies in
the Preservation of Cultural Heritage and the Tackling of Societal
Challenges''). T.R.\ also acknowledges the hospitality of (and partial support
from) FORTH in Heraklion, Crete.
}

\baselineskip=11pt

\bibliographystyle{plain}

\end{document}